\documentclass[a4paper,12pt]{article}
\begin{document}

\title{Structure and homomorphisms of microbundles.}
\author{Sergey Victor Ludkovsky.}
\date{16 March 2023}
\maketitle
\begin{abstract} This article is devoted to investigations of a structure and homomorphisms
of microbundles. Microbundles are generalizations of manifolds. For
manifolds it was studied when their families of homomorphism can be
supplied with the manifold structure. But for microbundles this
problem was not yet investigated. Continuous homomorphisms of
microbundles are studied. Topologies on families of homomorphisms of
microbundles are investigated. Relations of their structure with
microbundles are scrutinized. Necessary and sufficient conditions
are studied under which families of homomorphism of microbundles in
their turn can be supplied with a microbundle structure.
\footnote{key words and phrases: microbundle; structure; homomorphism; set-open topology  \\
Mathematics Subject Classification 2020: 55R60; 54B35; 54H13}

\end{abstract}
\par Address: Dep. Appl. Mathematics, MIREA - Russian Techn. University,
\par av. Vernadsky 78, Moscow 119454, Russia; e-mail: sludkowski@mail.ru

\section{Introduction.}
Topological manifolds and geometric bundles are very important in
mathematics and its applications (see, for example,
\cite{bourbalg,hewrossb,lhdtmta18,ludqimb}-\cite{luddfnafjms9,roo,sch1,xscscta19,weil}).
Their generalizations are microbundles. They are widely used in
topology and algebraic topology \cite{fedchigb,fedumn86,klmbim66,
milnor64,sulwintb,ludkmicrta19}. There were several works about
microbundles on ${\bf R}^n$ and Banach spaces over $\bf R$, but on
topological modules over rings they were little studied. Studies of
microbundles over topological rings are motivated by problems of
general topology, algebraic topology, noncommutative geometry,
representation theory, bundles over topological groups and group
rings \cite{avgmctrmb88,avmtrmb81,felldorb,hewrossb}, mathematical
physics. Microbundles have more complicated structure than
manifolds, or topological manifolds. For manifolds it was studied
when their families of homomorphism can be supplied with the
manifold structure. But for microbundles this problem was not yet
investigated.
\par We recall their definition and notation (or see \cite{ludkmicrta19}).
\par {\bf Notation I1.} Assume that $\bf F$ is a topological ring
such that its topology $\tau _{\bf F}$ is neither discrete nor
antidiscrete. In this case $({\bf F}, \tau _{\bf F})$ is called a
proper topological ring. A topological left module over $\bf F$ is
denoted by $X_{\bf F}$, or shortly by $X$ if $\bf F$ is specified.
The ring $\bf F$ is supposed to be associative and commutative
relative to the addition, but may be noncommutative or
nonassociative relative to the multiplication.
\par Henceforth, it will also be written a ring or a module
as a shortening of a topological ring or a topological module. Their
homomorphisms are supposed to be continuous. Neighborhoods in
topological spaces, modules, rings will be open if something other
will not be specified and the topological terminology is used in the
sense of the book \cite{eng}.

\par {\bf Definition I2.} Assume that the following conditions are
satisfied:
\par $(I2.1)$ there are are topological
spaces $A$ and $E$;
\par $(I2.2)$ there are continuous maps $i: A\to E$ and $p: E\to A$ such
that $p\circ i=id$, where $id=id_A: A\to A$ is the identity map,
$p\circ i$ denotes a composition of maps;
\par $(I2.3)$ for each $b\in A$ neighborhoods $U$ of $b$ and $V$ of $i(b)$ exist
such that $i(U)\subset V$ and $p(V)\subset U$ and $V$ is
homeomorphic to $U\times X$, where $h_V: V\to U\times X$ is a
homeomorphism; \par $(I2.4)$ there are continuous maps: a projection
$\hat{\pi }_1: U\times X\to U$, an injection $\iota _0: U\to U\times
X$, also a projection $\hat{\pi }_2: U\times X\to X$ such that
$\hat{\pi }_1(d,x)=d$, $ ~ \iota _0(d)=(d,0)$ and $\hat{\pi
}_2(d,x)=x$ for each $d\in U$ and $x\in X$; such that they are
supposed to satisfy the identity: $\hat{\pi }_1\circ \iota
_0|_U=p|_V \circ i|_U $, where $i|_U$ denotes the restriction of $i$
to $U$.
\par  It is said that the conditions $(I2.1)$-$(I2.4)$
define a microbundle ${\cal B}={\cal B}(A, E, {\bf F}, X, i, p)$
with a fibre $X=X_{\bf F}$ of ${\cal B}$, where $A$ is called a base
space, $E$ is called a total space.
\par If some data are specified, like $\bf F$
or $X$, they can be omitted from ${\cal B}(A, E, {\bf F}, X, i, p)$
in order to shorten the notation. Previously some their topological
properties were begun studied in \cite{ludkmicrta19}.
\par This article is devoted to investigations
of microbundle homomorphisms. Topologies on families of
homomorphisms of microbundles are investigated. Preliminaries on
them are described in Theorem 4, Proposition 6, Corollaries 7-10. An
existence of canonical microbundles is provided by Proposition 15.
Specific regularity of microbundles and families of their
homomorphisms over topological rings is investigated in Theorems 12,
22. A comparison of topologies is studied (see Theorem 23). Examples
24.$(1)$-$(5)$ provide some types of topologies. Relations of a
homomorphism family structure with microbundles are scrutinized (see
Theorem 28). Embeddings of microbundles are investigated in Theorem
30. In Theorems 32, 33, 34, 36 induced homomorphisms are studied.
Necessary and sufficient conditions are studied under which families
of homomorphism of microbundles in their turn can be supplied with a
microbundle structure.
\par All main results of this paper are obtained for the first time.
They can be used for subsequent investigations in topology,
algebraic topology, noncommutative geometry, representation theory,
mathematical analysis, mathematical physics and their applications.

\section{Homomorphisms of microbundles.}
\par {\bf Definition 1.} Let ${\cal B}_j={\cal
B}(A_j,E_j,F_j,X_j,i_j,p_j)$ be microbundles, $j\in {\bf N}$. If
$A_2\subset A_1$, $E_2\subset E_1$, $F_2\subset F_1$, $X_2\subset
X_1$, $i_1|_{A_2}=i_2$, $p_1|_{E_2}=p_2$, for each $b\in A_2$
neighborhoods $U_j$ of $b$ in $A_j$ and $V_j$ of $i_2(b)$ in $E_j$
exist such that conditions $(I2.3)$, $(I2.4)$  are satisfied,
$U_2=U_1\cap A_2$, $V_2=V_1\cap E_2$, $h_{V_1,1}|_{V_2}=h_{V_2,2}$,
$\iota _{0,1}|_{U_2}=\iota _{0,2}$, ${\hat \pi }_{1,1}|_{U_2\times
X_2}={\hat \pi }_{1,2}$, ${\hat \pi }_{2,1}|_{U_2\times X_2}={\hat
\pi }_{2,2}$, then ${\cal B}_2$ is called the submicrobundle in the
microbundle ${\cal B}_1$, where $h_{V_j,j}$, $\iota _{0,j}$,
$\hat{\pi }_{1,j}$, $\hat{\pi }_{2,j}$ correspond to ${\cal B}_j$.
\par If ${\cal B}_2$ is the submicrobundle in the microbundle ${\cal
B}_1$ and $A_2$ is open in $A_1$ put $\alpha _1=1$ otherwise $\alpha
_1=0$, if $E_2$ is open in $E_1$ put $\alpha _2=1$ otherwise $\alpha
_2=0$, if $F_2$ is open in $F_1$ put $\alpha _3=1$ otherwise $\alpha
_3=0$, if $X_2$ is open in $X_1$ put $\alpha _4=1$ otherwise $\alpha
_4=0$. Then it will be said that ${\cal B}_2$ is $\beta $-open in
${\cal B}_1$, where $\beta =(\alpha _1, 2\alpha _2, 3\alpha _3, 4
\alpha _4)$. Similarly is defined the $\beta $-closed ($\beta
$-compact) submicrobundle ${\cal B}_2$ in ${\cal B}_1$. For
shortening the notation zero coordinates of $\beta $ may be omitted.
For $\beta =(1,2,3,4)$ it may be said shortly the open (closed,
compact) submicrobundle instead of $\beta $-open ($\beta $-closed,
$\beta $-compact respectively).
\par Assume that $q=(b,e,f,x)$ with $b\in A_j$, $e\in E_j$, $f\in
F_j$, $x\in X_j$ for a fixed $j$ such that $i_j(b)=e_1\in E_j$,
$p_j(e_1)=b$, and conditions $(I2.3)$, $(I2.4)$ for neighborhoods
$U$ of $b$ in $A_j$ and $V$ of $e_1$ in $E_j$ are satisfied,
$\hat{\pi }_1\circ h_V(e)=b$, $\hat{\pi }_2\circ h_V(e)=x$. Then $q$
will be called an element in ${\cal B}_j$ and denoted by $q\in {\cal
B}_j$. It is written $C_j\subset {\cal B}_j$ if $(\forall q (q\in
C_j)\rightarrow (q\in {\cal B}_j))$ and $C_j$ is called a subset in
${\cal B}_j$. We also put $C_j=(C_j^1,C_j^2,C_j^3,C_j^4)$, where
$C_j^1\subset A_j$, $C_j^2\subset E_j$, $C_j^3\subset F_j$,
$C_j^4\subset X_j$. If in addition $C_j$ is a submicrobundle in
${\cal B}_j$ it will also be written $C_j\subset _M {\cal B}_j$,
that is there is an embedding $\pi $ from $C_j$ into ${\cal B}_j$
such that $\pi =id_{C_j}$. In particular, if $C_j = \{ q \} $, where
$q=(b,i_j(b),0,0)$, then it will also be written $q\in _M{\cal
B}_j$.
\par {\bf Remark 2.} If the subring $F_2$ is open in the ring $F_1$
(or the left $F_2$-module $X_2$ is open in the left $F_1$-module
$X_1$), then $F_2$ is closed in $F_1$ (or $X_2$ is closed in $X_1$
respectively) by Theorem 5.5 in \cite{hewrossb}, because $F_j$ (or
$X_j$ respectively) has the structure of the commutative topological
group relative to the addition. For topological spaces (topological
rings or left modules) $A$ and $B$ by $B^A={\sf C}(A,B)$ ($Hom(A,B)$
respectively) is denoted the family of all continuous maps
(continuous homomorphisms respectively) $f: A\to B$. If $A$ and $B$
are left $F$-modules, then $Hom_F(A,B)$ denotes the family of all
continuous left $F$-linear homomorphisms $g: A\to B$.
\par If the ring $F$ is associative, then it is as usually assumed
that $(ab)x=a(bx)$ for each $a$, $b$ in $F$, $x$ in the left
$F$-module $X=X_F$. If $F$ is associative and commutative, it also
is assumed that $a(bx)=b(ax)$ for each $a$, $b$ in $F$, $x$ in $X$.
\par By $Hom ({\cal B}_n,{\cal B}_k)$ is denoted a family of all
continuous homomorphisms $\pi ^n_k: {\cal B}_n \to {\cal B}_k$ of
the microbundles satisfying the conditions $(2.1)-(2.5)$:
\par $(2.1)$ $\pi ^n_k = (\pi ^{1,n}_{1,k}; \pi ^{2,n}_{2,k}; \pi
^{3,n}_{3,k}; \pi ^{4,n}_{4,k})$ with
\par $(2.2)$ $\pi ^{1,n}_{1,k}: A_n\to A_k$ and  $\pi ^{2,n}_{2,k}:
E_n\to E_k$, \par $\pi ^{3,n}_{3,k}: {\bf F}_n\to {\bf F}_k$ and
$\pi ^{4,n}_{4,k}: X_n\to X_k$, where
\par $(2.3)$ $\pi ^{2,n}_{2,k}\circ i_n=i_k\circ \pi
^{1,n}_{1,k}$ and $p_k\circ \pi ^{2,n}_{2,k}=\pi ^{1,n}_{1,k}\circ
p_n$, \par $(2.4)$ $h_k\circ \pi ^{2,n}_{2,k}\circ
h_n^{-1}(a_n,u_nx_n+v_ny_n)=$\par $ \pi ^{3,n}_{3,k}(u_n)h_k\circ
\pi ^{2,n}_{2,k}\circ h_n^{-1}(a_n,x_n) + \pi ^{3,n}_{3,k}(v_n)
h_k\circ \pi ^{2,n}_{2,k}\circ h_n^{-1}(a_n,y_n)$,
\par $(2.5)$ $\hat{\pi }_{2,k} \circ h_k\circ \pi ^{2,n}_{2,k}=\pi
^{4,n}_{4,k}\circ \hat{\pi }_{2,n} \circ h_n$ and \par $\pi
^{4,n}_{4,k}(u_nx_n+v_ny_n)=\pi ^{3,n}_{3,k}(u_n) \pi
^{4,n}_{4,k}(x_n)+ \pi ^{3,n}_{3,k}(v_n) \pi ^{4,n}_{4,k}(y_n)$ \\
for every $u_n$ and $v_n$ in ${\bf F}_n$, $x_n$ and $y_n$ in $X_n$,
$a_n\in A_n$, where $h_n: V_n\to U_n\times X_n$ is a local
homeomorphism for an open subset $V_n$ in $E_n$ corresponding to an
open neighborhood $U_n$ of a point $a_n$ in $A_n$, where $\pi
^{3,n}_{3,k}$ is the ring homomorphism, that is $\pi
^{3,n}_{3,k}(u_n+v_n)= \pi ^{3,n}_{3,k}(u_n)+ \pi ^{3,n}_{3,k}(v_n)$
and $\pi ^{3,n}_{3,k}(u_nv_n)= \pi ^{3,n}_{3,k}(u_n) \pi
^{3,n}_{3,k}(v_n)$ for every $u_n$ and $v_n$ in ${\bf F}_n$ (see
also Introduction or \cite{ludkmicrta19}). If $F_n=F_k=F$, then
\par $Hom_{F} ({\cal B}_n,{\cal B}_k) := \{ \pi ^n_k \in Hom ({\cal
B}_n,{\cal B}_k): \pi ^{4,n}_{4,k}\in Hom_{F}(X_n,X_k) \} $. \par If
${\cal B}={\cal B}(A,E,F,X,i,p)$ is the microbundle, $\tau _A$,
$\tau _E$, $\tau _F$, $\tau _X$ are topologies on the base space
$A$, on the total space $E$, on the topological ring $F$, on the
topological left $F$-module $X$ respectively such that conditions
$(I2.1)$-$(I2.4)$ are satisfied, then $\tau _{\cal B}=(\tau _A, \tau
_E, \tau _F, \tau _X)$ is called a topology on ${\cal B}$.
\par If $\pi ^2_1$ is bijective,
$\pi ^2_1: {\cal B}_2\to {\cal B}_1$ and $(\pi ^2_1)^{-1}: {\cal
B}_1\to {\cal B}_2$ are (continuous) homomorphisms of the
microbundles, then $\pi ^2_1$ is called an isomorphism of the
microbundles (particularly, for ${\cal B}_1={\cal B}_2$ an
automorphism).
\par Having $q\in C_l\subset {\cal B}_j$ for each $l\in \Lambda $,
by Definition 1, one naturally gets $C=\bigcup_{l\in \Lambda
}C_l\subset {\cal B}_j$ and $P=\bigcap_{l\in \Lambda }C_l\subset
{\cal B}_j$, where $\Lambda $ is a set.
\par {\bf Definition 3.} For the microbundles ${\cal B}_j$ assume
that $C_j\subset {\cal B}_j$, that is $C_j^1\subset A_j$,
$C_j^2\subset E_j$, $C_j^3\subset F_j$, $C_j^4\subset X_j$. We put
${\sf H}(C_j,C_l)={\sf H}_{{\cal B}_j,{\cal B}_l}(C_j,C_l)$ to be
the family of all $\pi ^j_l \in Hom ({\cal B}_j,{\cal B}_l)$ such
that $\pi ^{k,j}_{k,l}(C_j^k)\subset C_l^k$ for each $k\in \{ 1, 2,
3, 4 \} $, ${\sf H}_k(C_j,C_l):=\{ \pi ^{k,j}_{k,l}: \pi ^j_l\in
{\sf H}(C_j,C_l) \} $. Similarly if $F_j=F_l$, then ${\sf
H}_{F_j}(C_j,C_l)={\sf H}_{F_j;{\cal B}_j,{\cal B}_l}(C_j,C_l)
:=Hom_{F_j} ({\cal B}_j,{\cal B}_l)\cap {\sf H}(C_j,C_l) $. Then it
will be said that the collection of all families ${\sf H}(C_j,C_l)$
(or ${\sf H}_{F_j}(C_j,C_l)$ if $F_j=F_l$) with compact $C_j$ in
${\cal B}_j$ and open $C_l$ in ${\cal B}_l$ forms a prebase of a
compact-open topology $\tau _{co}= \tau _{co}({\cal B}_j,{\cal
B}_l)$ on $Hom ({\cal B}_j,{\cal B}_l)$ (or on $Hom_{F_j} ({\cal
B}_j,{\cal B}_l)$ respectively). If the microbundles ${\cal B}_j$
and ${\cal B}_l$ are specified, they may be omitted from ${\sf
H}_{{\cal B}_j,{\cal B}_l}(C_j,C_l)$ or ${\sf H}_{F_j;{\cal
B}_j,{\cal B}_l}(C_j,C_l)$ for shortening the notation.
\par {\bf Theorem 4.} {\it If ${\cal B}_3$ is a submicrobundle in a
$T_2$ microbundle ${\cal B}_1$, ${\cal B}_4$ is a closed
submicrobundle in a microbundle ${\cal B}_2$, $F_j$ and $X_j$ are
$T_0$, then ${\sf H}_{{\cal B}_1,{\cal B}_2}({\cal B}_3,{\cal B}_4)$
is closed in $(Hom({\cal B}_1,{\cal B}_2), \tau _{co})$. Moreover,
if $F_1=F_2=F_3=F_4$, then ${\sf H}_{F_1}({\cal B}_3,{\cal B}_4)$ is
closed in $(Hom_{F_1}({\cal B}_1,{\cal B}_2), \tau _{co})$.}
\par {\bf Proof.} Let $C_1$ and $C_3$ be compact in ${\cal B}_1$ and
${\cal B}_3$ respectively, $C_2$ and $C_4$ be open in ${\cal B}_2$
and ${\cal B}_4$ respectively. On the other hand, if $\pi ^3_4\in
{\sf H}({\cal B}_3,{\cal B}_4)$, then by Definitions 1 and 3 there
exists $\pi ^1_2\in Hom ({\cal B}_1,{\cal B}_2)$ such that $\pi
^1_2({\cal B}_3)\subset {\cal B}_4$ and $\pi ^1_2|_{{\cal B}_3}=\pi
^3_4$. In view of Theorems 4.8 in \cite{hewrossb} and 3.1.10 in
\cite{eng}, conditions $(2.1)$, $(2.2)$, $(2.5)$ above (or see
\cite{ludkmicrta19}) and Definition 1, $C_1^3C_1^4$ and $C_3^3C_3^4$
are compact in $X_1$ and $X_3$ respectively, $\pi
^{4,3}_{4,4}(C_3^4)$, $\pi ^{3,3}_{3,4}(C_3^3)\pi
^{4,3}_{4,4}(C_3^4)=\pi ^{4,3}_{4,4}(C_3^3C_3^4)$ and $\pi
^{4,1}_{4,2}(C_1^4)$, $\pi ^{3,1}_{3,2}(C_1^3)\pi
^{4,1}_{4,2}(C_1^4)=\pi ^{4,1}_{4,2}(C_1^3C_1^4)$ are compact in
$X_4$ and $X_2$ respectively, where $C^3_1C^4_1= \{ x\in X_1: ~
\exists u\in C^3_1, ~ \exists y\in C^4_1, ~ uy=x \} $. Then we
deduce that $\pi ^{3,1}_{3,2}(C_1^3)\subset C_2^3$, $\pi
^{4,1}_{4,2}(C_1^4)\subset C_2^4$, $\pi
^{4,1}_{4,2}(C_1^3C_1^4)\subset C_2^3C_2^4$, $\pi
^{3,3}_{3,4}(C_3^3)\subset C_4^3$, $\pi ^{4,3}_{4,4}(C_3^4)\subset
C_4^4$, $\pi ^{4,3}_{4,4}(C_3^3C_3^4)\subset C_4^3C_4^4$.
\par Note that $(Hom ({\cal B}_1,{\cal B}_2), \tau _{co})$ is the
closed topological subspace in the product topological space
$A_2^{A_1}\times E_2^{E_1}\times F_2^{F_1}\times X_2^{X_1}$ supplied
with the Tychonoff product topology, where $A_2^{A_1}$, $E_2^{E_1}$,
$F_2^{F_1}$, $X_2^{X_1}$ are supplied with the corresponding
compact-open topologies (see Remark 2 and Definition 3). Indeed,
conditions $(2.1)$-$(2.5)$ are for continuous maps $\pi
^{1,1}_{1,2}$, $\pi ^{2,1}_{2,2}$ and continuous homomorphisms $\pi
^{3,1}_{3,2}$, $\pi ^{4,1}_{4,2}$ with continuous $i_1, i_2, p_1,
p_2, h_1, h_2$, $\hat {\pi }_{2,1}$, $\hat {\pi }_{2,2}$. Then we
get that $A_2^{A_1}\setminus A_4^{A_3}= \bigcup \{ {\sf H}_1(\{ b \}
,A_2\setminus A_4): b\in A_3 \} $, $E_2^{E_1}\setminus E_4^{E_3}=
\bigcup \{ {\sf H}_2(\{ e \} ,E_2\setminus E_4): e\in E_3 \} $,
$F_2^{F_1}\setminus \tilde{{\sf H}}_3(F_3,F_4) = \bigcup \{
\tilde{{\sf H}}_3(\{ f \} ,F_2\setminus F_4): f\in F_3 \} $,
$X_2^{X_1}\setminus \tilde{{\sf H}}_4(X_3,X_4) = \bigcup \{
\tilde{{\sf H}}_4(\{ x \} ,X_2\setminus X_4): x\in X_3 \} $, where
${\sf H}={\sf H}_{{\cal B}_1,{\cal B}_2}$, ${\sf H}_k={\sf
H}_{k;{\cal B}_1,{\cal B}_2}$ for each $k$, $\tilde{{\sf
H}}_3(B_3,F_4)= \{ \xi \in Hom (F_1,F_2): \xi (B_3)\subset F_4 \} $
for each $B_3\subset F_3$, $\tilde{{\sf H}}_4(P_3,X_4)= \{ \eta \in
Hom (X_1,X_2): \eta (P_3)\subset X_4 \} $ for each $P_3\subset X_3$,
where $F_2^{F_1}=Hom (F_1,F_2)$, $X_2^{X_1}=Hom (X_1,X_2)$,  where
as usually $A\setminus B= \{ a\in A: a\notin B \} $ for subsets $A$,
$B$ of a set $\Omega $. If $q\in {\cal B}_3$, $q=(b,e,f,x)$, $b\in
A_3$, $e\in E_3$, $f\in F_3$, $x\in X_3$, then the singleton $ \{ q
\} $ is evidently compact (see Definition 1). Therefore, ${\sf
H}_1(\{ b \} ,A_2\setminus A_4)$ is open in $A_2^{A_1}$, ${\sf
H}_2(\{ e \} ,E_2\setminus E_4)$ is open in $E_2^{E_1}$,
$\tilde{{\sf H}}_3(\{ f \} ,F_2\setminus F_4)$ is open in
$F_2^{F_1}$, $\tilde{{\sf H}}_4(\{ x \} ,X_2\setminus X_4)$ is open
in $X_2^{X_1}$. Consequently, $A_2^{A_1}\setminus A_4^{A_3}$ is open
in $A_2^{A_1}$, $E_2^{E_1}\setminus E_4^{E_3}$ is open in
$E_2^{E_1}$, $F_2^{F_1}\setminus \tilde{{\sf H}}_3(F_3,F_4)$ is open
in $F_2^{F_1}$, $X_2^{X_1}\setminus \tilde{{\sf H}}_4(X_3,X_4)$ is
open in $X_2^{X_1}$ and hence $A_4^{A_3}\times E_4^{E_3}\times
\tilde{{\sf H}}_3(F_3,F_4) \times \tilde{{\sf H}}_4(X_3,X_4)$ is
closed in $A_2^{A_1}\times E_2^{E_1}\times F_2^{F_1} \times
X_2^{X_1}$. Therefore, ${\sf H}({\cal B}_3,{\cal B}_4)$ is closed in
$Hom ({\cal B}_1,{\cal B}_2)$, since ${\sf H}({\cal B}_3,{\cal
B}_4)=(A_4^{A_3}\times E_4^{E_3}\times \tilde{{\sf H}}_3(F_3,F_4)
\times \tilde{{\sf H}}_4(X_3,X_4))\cap Hom ({\cal B}_1,{\cal B}_2)$
and $Hom ({\cal B}_1,{\cal B}_2)$ is closed in $A_2^{A_1}\times
E_2^{E_1}\times F_2^{F_1} \times X_2^{X_1}$.
\par {\bf Definition 5.} A microbundle ${\cal B}={\cal
B}(A,E,F,X,i,p)$ is $T_k$ if and only if $A$, $E$, $F$, $X$ are
$T_k$ as topological spaces, where $k\in \{ 0, 1, 2, 3,
3\frac{1}{2}, 4, 5, 6 \} $.
\par {\bf Proposition 6.} {\it There exists a (forgetful) covariant functor
${\cal F}$ from the category ${\cal M}$ of microbundles into the
category $\cal T$ of topological spaces.}
\par {\bf Proof.} The category ${\cal M}$ has a family
of objects $Ob ({\cal M})$, where each object is a microbundle. For
two microbundles ${\cal B}_1$ and ${\cal B}_2$ in $Ob ({\cal M})$
there exists a family $Mor ({\cal B}_2,{\cal B}_1)=Hom ({\cal
B}_2,{\cal B}_1)$ of all homomorphisms of ${\cal B}_2$ into ${\cal
B}_1$ satisfying conditions $(2.1)$-$(2.5)$.
\par There exists a map $m: Mor ({\cal B}_3,{\cal B}_2)\times Mor ({\cal B}_2,{\cal
B}_1)\to Mor ({\cal B}_3,{\cal B}_1)$ for each ${\cal B}_j\in Ob
({\cal M})$, $j\in \{ 1, 2, 3 \} $, such that $m(\pi ^3_2,\pi
^2_1)=\pi ^2_1 \circ \pi ^3_2\in Mor ({\cal B}_3,{\cal B}_1)$ for
each $\pi ^3_2\in Mor ({\cal B}_3,{\cal B}_2)$ and $\pi ^2_1 \in Mor
({\cal B}_2,{\cal B}_1)$. For each ${\cal B}_j\in Ob ({\cal M})$
there exists a unit morphism $1_{{\cal B}_j}$ such that $1_{{\cal
B}_j}\in Mor ({\cal B}_j,{\cal B}_j)$, $1_{{\cal B}_1}\circ \pi
^2_1=\pi ^2_1\circ 1_{{\cal B}_2}=\pi ^2_1$ for each $\pi ^2_1 \in
Mor ({\cal B}_2,{\cal B}_1)$. Thus ${\cal M}$ is the category.
\par For each microbundle ${\cal B}\in Ob ({\cal M})$, ${\cal
B}={\cal B}(A,E,F,X,i,p)$, we consider a topological space $A\times
E\times F\times X=:P_{\cal B}$ with the Tychonoff product topology.
Let $\hat{\iota }: A\to P_{\cal B}$ be a map such that $\hat{\iota
}(b)=(b,i(b),0,\hat{\pi }_2(h_V(i(b)))$ for each $b\in A$, where $V$
is a neighborhood of $i(b)$ in $E$ satisfying condition $(I2.3)$.
Let also $\hat{p}: P_{\cal B}\to A$ be a map such that $p\circ
\theta _2=\hat{p}$, where $\theta _2: P_{\cal B}\to E$ is the
projection map, $\theta _2(b,e,f,x)=e$ for each $b\in A$, $e\in E$,
$f\in F$, $x\in X$. Therefore, from $(I2.2)$ it follows that
$\hat{p}\circ \hat{\iota }=id$, where $id=id_A: A\to A$ is the
identity map. Thus there exists a map ${\cal F}: Ob ({\cal M})\to Ob
({\cal T})$ such that ${\cal F}({\cal B}_j)=P_{{\cal B}_j}$, $j\in
{\bf N}$, where $\cal T$ is the category of topological spaces. \par
Consider any $\pi ^2_1\in Mor ({\cal B}_2,{\cal B}_1)$, $\pi ^3_2\in
Mor ({\cal B}_3,{\cal B}_2)$ for ${\cal B}_j\in Ob ({\cal M})$,
$j\in {\bf N}$. In view of $(2.1)$-$(2.5)$ there are continuous maps
$\pi ^{1,j+1}_{1,j}: A_{j+1}\to A_j$; $\pi ^{2,j+1}_{2,j}:
E_{j+1}\to E_j$; $\pi ^{3,j+1}_{3,j}: F_{j+1}\to F_j$; $\pi
^{4,j+1}_{4,j}: X_{j+1}\to X_j$, $j\in {\bf N}$, consequently, $\pi
^{1,j+1}_{1,j}\circ \pi ^{1,j+2}_{1,j+1}: A_{j+2}\to A_j$; $\pi
^{2,j+1}_{2,j}\circ \pi ^{2,j+2}_{2,j+1}: E_{j+2}\to E_j$; $\pi
^{3,j+1}_{3,j}\circ \pi ^{3,j+2}_{3,j+1}: F_{j+2}\to F_j$; $\pi
^{4,j+1}_{4,j}\circ \pi ^{4,j+2}_{4,j+1}: X_{j+2}\to X_j$.
Homomorphisms $\pi ^n_k: {\cal B}_n\to {\cal B}_k$ of microbundles
induce maps $\tilde{\pi }^n_k: P_{{\cal B}_n}\to P_{{\cal B}_k}$ by
$(2.1)$, $(2.2)$, where $\tilde{\pi }^n_k=\prod_{l=1}^4\pi
^{l,n}_{l,k}$ denotes a product map. From $(2.3)$ it follows that
$\tilde{\pi } ^{j+1}_j\circ \hat{\iota }_{j+1}=\hat{\iota }_j\circ
\pi ^{1,j+1}_{1,j}$ and $\hat{p}_j\circ \tilde{\pi }^{j+1}_j=\pi
^{1,j+1}_{1,j}\circ \hat{p}_{j+1}$. Condition $(2.5)$ implies that
$\hat{\pi }_{2,j}\circ h_j\circ \pi ^{2,j+1}_{2,j}\circ \theta
_{2,j+1}=\pi ^{4,j+1}_{4,j}\circ \hat{\pi }_{2,j+1}\circ
h_{j+1}\circ \theta _{2,j+1}$. Therefore, there exists a covariant
functor ${\cal F}$ such that ${\cal F}(\pi ^2_1\circ \pi ^3_2)={\cal
F}(\pi ^2_1)\circ {\cal F}(\pi ^3_2)$, ${\cal F}: {\cal M}\to {\cal
T}$, ${\cal F}(\pi ^{j+1}_j)\in Mor ({\cal F}({\cal B}_{j+1}),{\cal
F}({\cal B}_j))$. Particularly, ${\cal F}(id_{{\cal B}_1})=id_{{\cal
F}({\cal B}_1)}: P_{{\cal B}_1}\to P_{{\cal B}_1}$, where
$id_Y(y)=y$ for each $y\in Y$.
\par {\bf Corollary 7.} {\it Let ${\cal B}$ be a microbundle.
Then ${\cal B}$ is $T_i$, where $0\le i\le 3\frac{1}{2}$ if and only if
${\cal F}({\cal B})$ is $T_i$.}
\par {\bf Corollary 8.} {\it Assume that ${\cal B}$ is a microbundle
such that ${\cal F}({\cal B})$ is $T_i$, where $0\le i\le 6$. Then ${\cal B}$ is $T_i$.}
\par {\bf Proof.} The assertions of Corollaries 7 and 8 follow from
Theorem 2.3.11 in \cite{eng} and Proposition 6, Definition 5 above.
\par {\bf Corollary 9.} {\it A submicrobundle ${\cal B}_1$ in a
microbundle ${\cal B}_2$ is open (closed or compact) if and only if
${\cal F}({\cal B}_1)$ is open (closed or compact respectively) in ${\cal F}({\cal B}_2)$.}
\par {\bf Proof.} This follows from Proposition 2.3.1, Corollary
2.3.4 and Theorem 3.2.4 in \cite{eng}, Definition 1 and Proposition 6 above.
\par {\bf Corollary 10.} {\it If $\pi ^2_1: {\cal B}_2\to {\cal B}_1$ is an isomorphism
of microbundles ${\cal B}_2$ and ${\cal B}_1$, then ${\cal F}(\pi ^2_1):
{\cal F}({\cal B}_2)\to {\cal F}({\cal B}_1)$ is a homeomorphism of topological spaces.}
\par {\bf Remark 11.} There are homeomorphisms $h^3_f: F\to F$ and
$h^4_x: X\to X$ of the topological ring $F$ and the topological left
$F$-module $X$ considered as the topological spaces such that
$h^3_f(g)=f+g$, $h^4_x(y)=x+y$ for each $f$, $g$ in $F$, $x$, $y$ in
$X$. Notice that generally they are not homomorphisms in nontrivial
cases. Therefore, $h_{f,x}:=(id_A, id_E, h^3_f, h^4_x): A\times
E\times F\times X\to A\times E\times F\times X$ is the
homeomorphism. On the other hand, ${\cal B}(\{ b \} , \{ i(b) \}
\cup \{ e \} , 0,0,i,p)$ is the trivial submicrobundle in the
microbundle ${\cal B}(A,E,F,X,i,p)$ for each $b\in A$ and $e\in E$.
In view of this and Corollary 10 we put ${\cal F}_1( \{ q \}
)=h_{q_3,q_4}({\cal F} ({\cal B}( \{ q_1 \} , \{ i(q_1) \} \cup \{
q_2 \} , 0,0,i,p))\setminus  h_{q_3,q_4} ({\cal F}({\cal B}( \{ q_1
\} , \{ i(q_1 ) \} , 0, 0,i,p))$ for each $q\in {\cal B}={\cal
B}(A,E,F,X,i,p)$, where $q=(q_1,q_2,q_3,q_4)$, $q_1\in A$, $q_2\in
E$, $q_3\in F$, $q_4\in X $ and conditions $(I2.1)$-$(I2.4)$ are
satisfied (see also Definition 1).
\par Using this we put
${\cal B}(\{ b \} , \{ e _1 \} , F, \{ 0 \} ,i_y,p_y)=:L_y{\cal B}(
\{ b \} , \{ e \} , F , \{ 0 \} , i,p)$ if and only if $b\in A$,
$e\in E$, $e_1\in E$, $p(e)=b$, $i(b)=e$, $p_y(e_1)=b$,
$i_y(b)=e_1$, $\hat{\pi }_2(h_V(e_1))=y\in X$, $\hat{\pi
}_2(h_V(e))=0\in X$, where $\hat{\pi }_2(h_V\circ i_y)=h^4_y\circ
\hat{\pi }_2(h_V\circ i)$, $p_y\circ i_y=id$, $e\in V$ and
conditions $(I2.1)$-$(I2.4)$ are satisfied, ${\cal B}( \{ b \} , \{
e \} , F , \{ 0 \} , i,p)\subset _M{\cal B}(A,E,F,X,i,p)$. This
implies that $L_y\circ L_{-y}=L_{-y}\circ L_y=id_{\cal B}.$ Thus an
embedding $\hat{L}_y: {\cal B}(\{ b \} , \{ e _1 \} , F, \{ 0 \}
,i_y,p_y)\hookrightarrow {\cal B}(A,E,F,X,i,p)$ exists induced by a
shift operator $L_{-y}$.

\par {\bf Theorem 12.} {\it If the microbundle ${\cal B}_1$ is
$T_1\cap T_3$, then $(Hom ({\cal B}_2,{\cal B}_1), \tau _{co})$ is
$T_1\cap T_3$.}
\par {\bf Proof.} Take any $\pi ^2_1\in {\sf H}(C_2,C_1)$, where $C_2$ is
compact in ${\cal B}_2$ and $C_1$ is open in ${\cal B}_1$. In view
of Proposition 6 and Corollary 9 there exists an open subset $V$ in
${\cal B}_1$ such that $\pi ^2_1\in {\sf H}(C_2,V)$ and $cl_Y[{\sf
H}(C_2,V)]\subset {\sf H}(C_2,C_1)$, where $cl_Y[{\sf H}]$ denotes
the closure of ${\sf H}$ in $Y=(Hom ({\cal B}_2,{\cal B}_1), \tau
_{co})$; $cl_{{\cal B}_1}V\subset C_1$, since $\pi ^2_1(C_2)$ is
compact by Theorem 3.1.10 in \cite{eng}.  We have $\pi
^2_1(C_2)\subset C_1$ and ${\cal B}_1$ is $T_1\cap T_3$.
\par To each $q\in {\cal B}_2$ there corresponds
$q=(q_1,q_2,q_3,q_4)$ with $q_1\in A_2$, $q_2\in E_2$, $q_3\in F_2$,
$q_4\in X_2$, consequently, ${\cal F}( \{ q \} )\in {\cal F}({\cal
B}_2)$ by Proposition 6 and Remark 11. Therefore, for each $q\in
C_2$ there exists an open set $W_q$ in ${\cal B}_1$ such that $\pi
^2_1(q)\in W_q$ and $cl_{{\cal B}_1} W_q\subset C_1$ by Corollary 9.
The open covering $ \{ W_q: q\in C_2 \} $ of $\pi ^2_1(C_2)$
contains a finite subcovering $ \{ W_{q_1},...,W_{q_n} \} $, where
$n\in {\bf N}$. We put $W=\bigcup_{k=1}^n W_{q_k}$, hence $\pi
^2_1(C_2)\subset W$ and $cl_{{\cal B}_1}W=\bigcup_{k=1}^ncl_{{\cal
B}_1}W_{q_k}\subset C_1$. This implies that $\pi ^2_1\in {\sf
H}(C_2,W)$ and ${\sf H}(C_2,cl_{{\cal B}_1}W)\subset {\sf
H}(C_2,C_1)$. By virtue of Theorem 4  ${\sf H}(C_2,cl_{{\cal
B}_1}W)$ is closed, consequently, $cl_{Hom ({\cal B}_2,{\cal
B}_1)}{\sf H}(C_2,W)\subset cl_{Hom ({\cal B}_2,{\cal B}_1)} {\sf
H}(C_2,cl_{{\cal B}_1}W)= {\sf H}(C_2,cl_{{\cal B}_1}W)\subset {\sf
H}(C_2,C_1)$. This implies that $(Hom ({\cal B}_2,{\cal B}_1), \tau
_{co})$ is $T_1\cap T_3$.
\par {\bf Corollary 13.} {\it Let ${\cal B}_j$ be a microbundle for
each $j\in \{ 1, 2, 3 \} $, let $\pi ^2_1\in Hom ({\cal B}_2,{\cal B}_1)$
and $\eta (\pi ^3_2):=\pi ^2_1\circ \pi ^3_2$  for each $\pi ^3_2\in Hom ({\cal B}_3,{\cal B}_2)$.
Then $\eta :   Hom ({\cal B}_3,{\cal B}_2)\to  Hom ({\cal B}_3,{\cal B}_1)$ is a continuous map
relative to the compact open-topology on $Hom ({\cal B}_3,{\cal B}_2)$ and $Hom ({\cal B}_3,{\cal B}_1)$.}
\par {\bf Proof.} If $C_j\subset {\cal B}_j$ for each $j\in \{ 1, 2,
3 \} $, then the inclusion $\pi ^2_1(\pi ^3_2(C_3))\subset C_1$ is
equivalent to $\pi ^3_2(C_3)\subset (\pi ^2_1)^{-1}(C_1)$, since
$\pi ^{k,2}_{k,1}(\pi ^{k,3}_{k,2}(C^k_3))\subset C_1^k$ is
equivalent to $\pi ^{k,3}_{k,2}(C_3^k)\subset (\pi
^{k,2}_{k,1})^{-1}(C_1^k)$ for each $k\in \{ 1, 2, 3, 4 \} $.
Therefore, \par $\eta ^{-1}({\sf H}_{{\cal B}_3,{\cal
B}_1}(C_3,C_1))={\sf H}_{{\cal B}_3,{\cal B}_1}(C_3,(\pi
^2_1)^{-1}(C_1))$.
\par {\bf Definition 14.} Assume that a microbundle ${\cal B}={\cal
B}(A,E,F,X,i,p)$ satisfies the condition: for each $e\in E$ there
exists an open neighborhood $V_e$ of $e$ in $E$ and an open
neighborhood $U_b$ of $b=p(e)$ in $A$ such that $V_e$ is
homeomorphic to $U\times X$, where conditions $(I2.1)$-$(I2.4)$ are
satisfied. Then ${\cal B}$ is called a canonical microbundle.
\par {\bf Proposition 15.} {\it If ${\cal B}_2={\cal
B}(A,E_2,F,X,i_2,p_2)$ is a microbundle, then there exists a canonical microbundle
${\cal B}_1={\cal B}(A,E_1,F,X,i_1,p_1)$ base neighborhood isomorphic with ${\cal B}_2$
and ${\cal B}_1\subset _M{\cal B}_2$.}
\par {\bf Proof.} We put \par $(15.1)$ $E_1=\bigcup \{ e\in E_2: \exists b\in A,
\exists U\in \tau _A, b\in U, \exists V\in \tau _{E_2}, e\in V,
\exists h_V: V\to U\times X \mbox{ is a homeomorphism }, \& (I2.3)
\& (I2.4) \mbox{ are satisfied } \} $,
\\ where $\tau _A$ and $\tau _{E_2}$ denote topologies on $A$ and $E_2$
respectively.   Therefore, $E_1\subset E_2$ and $i_2(A)\subset E_1$,
so we put $i_2=i_1$. By the definition of the microbundle $p_2\circ
i_2=id_A$, $p_2(E_2)=A$. Supply $E_1$ with the topology $\tau
_{E_1}$ inherited from $E_2$ and put $p_1=p_2|_{E_1}$. Therefore,
$p_1\circ i_1=id_A$, $id_A: A\to A$, $p_1(E_1)=A$, $p_1|_V\circ
i_1|_U=\hat{\pi }_1\circ \iota _0|_U$ for $U$ and $V$ as above in
$(15.1)$.   Thus conditions $(I2.1)$-$(I2.4)$ are satisfied and
${\cal B}_1$ is the microbundle. Since $A_2=A_1=A$, $X_2=X_1=X$,
$F_2=F_1=F$, $E_1\subset E_2$, then ${\cal B}_1$ is the
submicrobundle in ${\cal B}_2$, ${\cal B}_1\subset _M{\cal B}_2$.
With $V_1=E_1=V_2$, $g=id_{V_1}$ and $s=id_A$, one gets that ${\cal
B}_1$ and ${\cal B}_2$ are base neighborhood isomorphic by
Definition 4 in \cite{ludkmicrta19}.
\par {\bf Definition 16.} Let $X_1$ and $X_2$
be topological left $F$-modules, where $F$ is a topological ring.
The family of all continuous left $F$-linear homomorphisms $f:
X_2\to X_1$ is denoted by $Hom_F(X_2,X_1)$.
\par Let $q: X_{2,F_2} \to  X_{1,F_1}$ be a
continuous homomorphism satisfying $(16.1)$-$(16.3)$:
\par  $(16.1)$ $q(x+y)=q(x)+q(y)$,
\par $(16.2)$ $q(ax)=g(a)q(x)$
\\ for each $a\in F_2$, $x$ and $y$ in $X_2$, where $g=g_q: F_2\to
F_1$ is a continuous homomorphism
\par $(16.3)$ $g(a+b)=g(a)+g(b)$, $g(ab)=g(a)g(b)$ \\ for each $a$ and
$b$ in $F_2$, where $X_j=X_{j,F_j}$ is the topological left
$F_j$-module, $F_j$ is the topological ring $j\in \{ 1, 2 \} $. The
family of all continuous homomorphisms satisfying conditions
$(16.1)$-$(16.3)$ will be denoted by $Hom (X_2, X_1)$.
\par We recall, that if $X$ is a topological space, $\delta $ is a
nonempty family of nonvoid subsets of $X$ such that for each $x\in
X$ and an open neighborhood $U$ of $x$ there exists $B\in \delta $
with $x\in B\subset U$, then $\delta $ is called a network on $X$.
The network $\delta $ is called closed (or compact), if each $B\in
\delta $ is closed (or compact respectively) in $X$.
\par Let $\delta _j$ be a closed (or a compact network) on $X_j$ and
$ \alpha _j$ be a closed (or a compact network respectively) on
$F_j$. Then $Hom^{\alpha _2, \delta _2}(X_{2,F_2}, X_{1,F_1})$ is
$Hom (X_{2,F_2}, X_{1,F_1})$ supplied with the set-open (or
compact-open respectively) topology denoted by  $\tau
(F_2,F_1,X_2,X_1,\alpha _2,\delta _2)$. Similarly $Hom_F^{\alpha _2,
\delta _2}(X_{2,F},X_{1,F})$ is $Hom _F(X_{2,F}, X_{1,F})$ supplied
with the set-open (or compact-open respectively) topology $\tau
(F,X_2,X_1,\alpha _2,\delta _2)$. If some data are specified, they
may be omitted in this notation for its shortening.
\par Let a topological ring $F$ be fixed. Henceforth canonical microbundles are considered if some other
will not be specified. Let ${\cal B}={\cal B}(A,E,F,X_F,i,p)$ be the
microbundle. Let $\alpha $ be a family of subsets in ${\cal B}$ (in
particular they may be submicrobundles over the same ring $F$). A
nonempty family $\delta $ of nonvoid subsets in ${\cal B}$ is called
an $\alpha $-$M$-network on ${\cal B}$ provided that for each $C\in
\alpha $ and an open subset ${\cal U}$ in ${\cal B}$ there exists
${\cal Q}\in \delta $ such that $C\subset {\cal Q}\subset {\cal U}$.
\par A $M$-network on ${\cal B}$ is called an $\alpha _0$-$M$-network on
${\cal B}$, if it satisfies the following condition: $((C\in \alpha
_0)\leftrightarrow (\exists q\in {\cal B}, C=\{ q \} ))$ (see
Definition 1). The $\alpha $-$M$-network is called $\gamma $-closed
(or $\gamma $-compact) provided that each member is $\gamma $-closed
(or $\gamma $-compact respectively).
\par If $C_k\subset {\cal B}$, ${\cal A}_k\subset _M{\cal B}$,
$D_j\subset {\cal R}$, ${\cal B}_j\subset _M{\cal R}$, where ${\cal
A}_k$, ${\cal B}$, ${\cal B}_j$, ${\cal R}$ are microbundles over
the same topological ring $F$, ${\cal B}={\cal B}(A,E,F,X_F,i,p)$,
${\cal R}={\cal B}(A',E',F,X'_F,i',p')$, ${\cal B}_j={\cal
B}(A_j',E_j',F,X'_{j,F},i_j',p_j')$, then for ${\cal A}_k={\cal
B}(A_k,E_k,F,Y_k,i,p)$ we put \par $(16.4)$ $ {\cal A}_1\biguplus
{\cal A}_2 ={\cal B}(A_1\cup A_2, E_1\biguplus E_2, F, Y, i,p)$, \\
where $A_k\subset A$, $E_k\subset E$, $Y_k\subset X$ for each $k$,
$Y=cl_X(Y_1+Y_2)$, $Y_k=Y_{k,F}$, $Y=Y_F$;
\par $(16.5)$ $E_1\biguplus E_2= \bigcup \{ h_V^{-1}((U_1\cup U_2)\times Y): b\in A_1\cup A_2, U_j\subset A_j, U_j\subset U,
i(b)\in V, j\in \{ 1, 2 \} \} $, where conditions $(I2.1)$-$(I2.4)$
are satisfied. In particular, ${\sf H}(L_y{\cal B}( \{ b \} , \{ e
\} , F , \{ 0 \} , i,p)), {\cal B}_j)$ may be also shortly denoted
by ${\sf H}( (b,y), {\cal B}_j)$.

\par {\bf Lemma 17.} {\it Let ${\cal A}_k$ and ${\cal B}_j$ be as in Definition 16.
Then \par $(17.1)$ ${\sf H}({\cal A}_k,{\cal B}_1\cap {\cal B}_2)={\sf H}({\cal A}_k,{\cal B}_1)
\cap {\sf H}({\cal A}_k,{\cal B}_2)$;
\par $(17.2)$ ${\sf H}_F({\cal A}_k,{\cal B}_1\cap {\cal B}_2)={\sf H}_F({\cal A}_k,{\cal B}_1)\cap {\sf H}_F({\cal A}_k,{\cal B}_2)$;

\par moreover, if $X'$ is complete as the topological commutative group, then
\par $(17.3)$ ${\sf H}({\cal A}_1\biguplus {\cal A}_2,{\cal B}_j)={\sf H}({\cal A}_1,{\cal
B}_j)\cap {\sf H}({\cal A}_2,{\cal B}_j)$ and
\par $(17.4)$ ${\sf H}_F({\cal A}_1\biguplus {\cal A}_2,{\cal B}_j)={\sf H}_F({\cal A}_1,{\cal B}_j)
\cap {\sf H}_F({\cal A}_2,{\cal  B}_j)$.}
\par {\bf Proof.} We consider
${\cal R}={\cal B}(A',E',F,X',i',p')$  and $Y_{1,2}'=X'_{1,F}\cap
X'_{2,F}$, where $X'_{j,F}\subset X'_F$ and ${\cal B}_j={\cal
B}(A_j',E_j',F,X_j',i_j',p_j')\subset _M{\cal R}$ for each $j$,
hence $FY_{1,2}'\subset Y_{1,2}'$ and $Y_{1,2}'$ is a commutative
group, consequently, $Y_{1,2}'$ is a topological left $F$-module. We
put $A_{1,2}'=A_1'\cap A_2'$ and $E_{1,2}'=E_1'\cap E_2'$, hence
$p_j'\circ i_j'|_{A_{1,2}'}= id|_{A_{1,2}'}$,
$i_j'|_{A_j'}=i'|_{A_j'}$, $p_j'|_{E_j'}=p'|_{E_j'}$ for each $j\in
\{ 1, 2 \} $, consequently, $p_1'|_{E_{1,2}'}= p_2'|_{E_{1,2}'}$,
$i_1'|_{A_{1,2}'}= i_2'|_{A_{1,2}'}$. If $b_j\in A_j'$,
$i_j'(b_j')=e_j'\in E_j'$, $e_j'\in V_j'$, $h_{V_j'}: V_j'\to
U_j'\times X_j'$ is a homeomorphism, where $U_j'$ is open in $A_j'$
for each $j$, then $h_{V'}|_{V_j'}=h_{V_j'}$, $V_j'\subset V'$,
$V_{1,2}'=V_1'\cap V_2'$, $h_{V'}|_{V_{1,2}'}: V_{1,2}'\to
U_{1,2}'\times Y_{1,2}'$ is the homeomorphism and conditions
$(I2.1)$-$(I2.4)$ are satisfied, since $(U_1'\times X_1')\cap
(U_2'\times X_2')=U_{1,2}'\times Y_{1,2}'$, where $U_{1,2}'=U_1'\cap
U_2'$. Therefore ${\cal B}_1\cap {\cal B}_2= {\cal
B}(A_{1,2}',E_{1,2}',F,Y_{1,2}',i',p')\subset _M{\cal R}={\cal
B}(A',E',F,X',i',p')$. This implies that \par ${\sf H}({\cal
A}_k,{\cal B}_1\cap {\cal B}_2)={\sf H}({\cal A}_k,{\cal B}_1)\cap
{\sf H}({\cal A}_k,{\cal B}_2)$ and
\par ${\sf H}_F({\cal A}_k,{\cal B}_1\cap {\cal B}_2)={\sf H}_F({\cal A}_k,{\cal B}_1)
\cap {\sf H}_F({\cal A}_k,{\cal B}_2)$.
\par We have
 $ {\cal A}_1\biguplus {\cal A}_2 \subset _M{\cal B}$. From Definition 26 and
 Theorem 27 in \cite{ludkmicrta19} it follows that $f({\cal A}_k)\subset _M{\cal B}_j$. If $\pi \in
 Hom({\cal A}_k,{\cal B}_j)$ for each $k\in \{ 1, 2 \} $, then $\pi
 ^{1}_{1}(A_1\cup A_2)=\pi ^1_1(A_1)\cup \pi ^1_1(A_2)\subset A_j'$, where
 ${\cal B}_j={\cal B}(A_j',E_j',F,X_j',i_j',p_j')\subset _M{\cal R}$. On
 the other hand,
 $\pi ^3_3(F)\subset F$, $\pi ^4_4(Y_1)\subset X'$ and $\pi ^4_4(Y_2)\subset X'$, consequently,
 by continuous extension $\pi ^4_4(Y)\subset X'$ \cite{eng,hewrossb}, since $\pi ^4_4\in Hom (Y_j,X')$ for each $j\in \{ 1, 2 \} $ and
$X'$ is complete as the topological commutative group, where
$Y=cl_X(X_1+X_2)$. Therefore, $\pi ^2_2(E_1\biguplus E_2)\subset
E'$, consequently, $\pi ({\cal A}_1\biguplus {\cal A}_2) \subset
{\cal R}$. Thus $\pi \in {\sf H}({\cal A}_1,{\cal B}_j)\cap {\sf
H}({\cal A}_2,{\cal B}_j)$ implies that $\pi \in {\sf H}({\cal
A}_1\biguplus {\cal A}_2,{\cal B}_j)$. \par Vice versa, if $\pi \in
{\sf H}({\cal A}_1\biguplus {\cal A}_2,{\cal B}_j)$, then $\pi
|_{{\cal A}_k}\in {\sf H}({\cal A}_k,{\cal B}_j)$ for each $k\in \{
1, 2 \} $, since ${\cal A}_k\subset _M{\cal A}_1\biguplus {\cal
A}_2$. Thus
\par ${\sf H}({\cal A}_1\biguplus {\cal A}_2,{\cal B}_j)={\sf H}({\cal A}_1,{\cal
B}_j)\cap {\sf H}({\cal A}_2,{\cal B}_j)$ and similarly
\par ${\sf H}_F({\cal A}_1\biguplus {\cal A}_2,{\cal B}_j)={\sf H}_F({\cal A}_1,{\cal B}_j)
\cap {\sf H}_F({\cal A}_2,{\cal B}_j)$.

\par {\bf Definition 18.} A topology on $Hom ({\cal B},{\cal R})$ is
called a $M$-set-$(\beta , \gamma )$-open topology provided that
there exists a $\beta $-closed $M$-network $\alpha $ on ${\cal B}$
such that $ \{ {\sf H}(C_k, D_j): C_k\in \alpha , D_j \mbox{ is a }
\gamma-\mbox{open subset in } {\cal R} \} $ is a subbase for the
topology $\tau _{\alpha , \beta , \gamma }$. Then $Hom ({\cal
B},{\cal R})$ supplied with the topology $\tau _{\alpha , \beta ,
\gamma }$ will be denoted by $Hom_{\alpha , \beta , \gamma } ({\cal
B},{\cal R})$. If $\beta =(1,2,3,4)$ or $\gamma =(1,2,3,4)$ it may
be omitted for shortening the notation.

\par {\bf Remark 19.} If ${\cal B}={\cal B}(A,E,F,X_F,i,p)$
and ${\cal R}={\cal B}(A',E',F',X'_{F'},i',p')$ are microbundles,
$B\subset {\cal B}$, $C\subset {\cal B}$, $D\subset {\cal R}$,
$G\subset {\cal R}$, then from Definition 3 it follows, that
\par $(19.1)$ ${\sf H}_{{\cal B},{\cal R}}(C,D\cap G)={\sf H}_{{\cal B},{\cal R}}(C,D)\cap
{\sf H}_{{\cal B},{\cal R}}(C,G)$,
\par $(19.2)$ ${\sf H}_{{\cal B},{\cal R}}(B\cup C,D)={\sf H}_{{\cal B},{\cal R}}(B,D)\cap
{\sf H}_{{\cal B},{\cal R}}(C,D)$.
\par Moreover, if for the microbundles ${\cal B}$ and ${\cal R}$ a
ring is the same, $F=F'$, then
\par $(19.3)$ ${\sf H}_{F;{\cal B},{\cal R}}(C,D\cap G)={\sf H}_{F;{\cal B},{\cal R}}(C,D)\cap
{\sf H}_{F;{\cal B},{\cal R}}(C,G)$,
\par $(19.4)$ ${\sf H}_{F;{\cal B},{\cal R}}(B\cup C,D)={\sf H}_{F;{\cal B},{\cal R}}(B,D)\cap
{\sf H}_{F;{\cal B},{\cal R}}(C,D)$.

\par {\bf Corollary 20.} {\it Let ${\cal B}_1={\cal B}(A_1,E_1,F_1,X_1,i_1,p_1)$
be a microbundle, where $X_1=X_{1,F_1}$. Let $F_2$ be a topological
ring and $X_2=X_{2,F_2}$ be a topological left $F_2$-module and
$f=(\pi ^{3,1}_{3,2}, \pi ^{4,1}_{4,2}): X_1\to X_2$ be an injective
homomorphism such that $\pi ^{3,1}_{3,2}: F_1\to F_2$ is a ring
homomorphism, $\pi ^{4,1}_{4,2}: X_1\to X_2$, $\pi
^{4,1}_{4,2}(u_1x_1+v_1y_1)=\pi ^{3,1}_{3,2}(u_1)\pi
^{4,1}_{4,2}(x_1)+ \pi ^{3,1}_{3,2}(v_1)\pi ^{4,1}_{4,2}(y_1)$ for
each $u_1$, $v_1$ in $F_1$, $x_1$, $y_1$ in $X_1$. Then there exists
a microbundle ${\cal B}_2={\cal B}(A_2,E_2,F_2,X_2,i_2,p_2)$ and an
injective homomorphism $\pi ^1_2: {\cal B}_1\to {\cal B}_2$, where
$A_2=A_1$. }
\par {\bf Proof.} In view of Theorem 32 in \cite{ludkmicrta19} the microbundle
${\cal B}_2={\cal B}(A_2,E_2,F_2,X_2,i_2,p_2)$ and the homomorphism
$\pi ^1_2: {\cal B}_1\to {\cal B}_2$ exist  satisfying conditions
$(2.1)$-$(2.5)$ such that $\pi^{1,1}_{1,2}=id_{A_1}$, where
$A_2=A_1$. Evidently, if $f$ is injective, then $\pi ^1_2$ is
injective, since $\pi^{1,1}_{1,2}=id_{A_1}$ and $\pi^{2,1}_{2,2}$ is
injective by $(32.2)$-$(32.5)$ in \cite{ludkmicrta19} and $(2.3)$,
$(2.4)$.

\par {\bf Definition 21.} Let $X$ be a topological space, let $F$ be a topological
ring, $r\in F$, $r\ne 0$. If for each $q\ne s$ in $X$ there exists a
continuous map $\mu : F\to X$ such that $\mu (0)=s$, $\mu (r)=q$,
then $X$ will be called $(F,r)$-arcwise connected.
\par A topological left $F$-module $X=X_F$ will be called completely
$(F,r)$-regular, if and only if for each $x\ne 0$ in $X$ there
exists $\psi \in Hom_F(X,X)$ such that $\psi (x)=r$, $X=(Fx)+Ker
(\psi )$, $(Fx)\cap Ker (\psi )= \{ 0 \} $, where $r\ne 0$, $r\in
F$, $r$ is fixed, $F$ and $X$ are nontrivial. The microbundle ${\cal
B}={\cal B}(A,E,F,X,i,p)$ will be called completely $(F,r)$-regular,
if and only if ${\cal B}$ is $T_1$ and for each closed subset $G$ in
${\cal B}$ and for each $q\in {\cal B}$, $q\notin G$, there exists
$\phi _{G,q}=\phi \in Hom ({\cal B},{\cal F})$ such that $\phi (G)=
\{ (0,0,0,0) \} $, $\phi (q)= (0,r,q_3,r) $, and $X$ is completely
$(F,r)$-regular, where ${\cal F}={\cal B}(0,F,F,F,i_F,p_F)$ with
$i_F(v)=v$ and $p_F(v)=v$ for each $v\in F$, where $r\ne 0$, $r\in
F$, $~ 0\in F$, $ ~ r$ is fixed.
\par The microbundle ${\cal R}$ will be called $({\cal
F},r)$-arcwise connected, if for each $q\ne s$ with $s_3=0$ and
$s_4=0$ in ${\cal R}$ there exists $\mu _{q,s}=\mu \in Hom ({\cal
F},{\cal R})$ such that $\mu (0,0,0,0)=s$ and $\mu (0,r,q_3,r)=q$,
where $r\ne 0$, $r\in F$, $r$ is fixed, $s=(s_1,s_2,s_3,s_4)$ (see
Remark 11). If $F$ is unital and $r=1$, then it will also be written
shortly $F$-arcwise connected and $F$-regular instead of
$(F,r)$-arcwise connected and $(F,r)$-regular correspondingly.

\par {\bf Theorem 22.} {\it The microbundle ${\cal
B}_1={\cal B}(A,E,F,X,i,p)$ is completely $(F,r)$-regular, if and
only if ${\cal B}_2={\cal B}(A,A,F,0,i_2,p_2)$ and $X$ are
completely $(F,r)$-regular, where $i_2=id_A$, $p_2=id_A$, $X=X_F$.}
\par {\bf Proof.} From the conditions of this theorem it follows
that ${\cal B}_2$ is isomorphic to a submicrobundle ${\cal B}_3=
{\cal B}(A,i(A),F,0,i,p|_{i(A)})$ in ${\cal B}_1$, ${\cal
B}_3\subset _M{\cal B}_1$, there is an isomorphism $\pi ^2_3: {\cal
B}_2\to {\cal B}_3$, since $p\circ i=id_A: A\to A$. For each $b\in
A$ there exists an open neighborhood $U_b$ of $b$ in $A$ and an open
neighborhood $V_b$ of $i(b)$ in $E$ and a homeomorphism $h_{V_b}:
V_b\to U_b\times X$, where $X=X_F$, by $(I2.3)$. Therefore,
$V_b\setminus \{ i(b) \} $ is open in $E$, since ${\cal B}_1$ is
$T_1$. Hence ${\cal B}_3$ is closed in ${\cal B}_1$ such that
$i_3=i$, $p_3=p|_{i(A)}$, $E_3=i(A)$.
\par If ${\cal B}_1$ is completely $(F,r)$-regular, then from
Definition 21, one gets that ${\cal B}_3$ is completely
$(F,r)$-regular, consequently, ${\cal B}_2$ is completely
$(F,r)$-regular.
\par Vice versa assume that ${\cal B}_2$ and $X$ are completely
$(F,r)$-regular, ${\cal G}={\cal B}(A_4,E_4,F,Y_F,i_4,p_4)$ is a
closed submicrobundle in ${\cal B}_1$, $q\in {\cal B}_1$, $q\notin
{\cal G}$, where $i_4=i|_{A_4}$, $p=p|_{E_4}$, $A_4\hookrightarrow
A$, $E_4\hookrightarrow E$, $Y_F\hookrightarrow X_F$.
\par In view of Theorem 32 in
 \cite{ludkmicrta19} and Corollary 20 there exists a continuous surjective
homomorphism $\pi ^1_3: {\cal B}_1\to {\cal B}_3$ induced by
$f=(id_F, \pi ^{4,1}_{4,3})$ such that $\pi ^{1,1}_{1,3}=id_A$,
where $\pi ^{4,1}_{4,3}(x)=0$ for each $x\in X$. If $A_4\ne A$,
$q_1\notin A_4$, where $q=(q_1,q_2,q_3,q_4)$, then from the complete
$(F,r)$-regularity of ${\cal B}_3$ using $\pi ^1_3$ we infer that
there exists $\phi \in Hom ({\cal B}_1,{\cal F})$ such that $\phi
({\cal G})= \{ (0,0,0,0) \} $ and $\phi (q)=(0,r,q_3,r)$. If $q_1\in
A_4$, then there exists a homeomorphism $h_{V_{q_1}}(q_2)=q_4$ such
that $q_4$ is nonzero $q_4\in X\setminus Y $ according to $(I2.3)$,
$(I2.4)$, since $q\notin {\cal G}$.
 By Definition 21 there exists
$\psi \in Hom_F(X,X)$ such that $\psi (q_4)\ne 0$, $X=(Fq_4)+Ker
(\psi )$ and $(Fq_4)\cap Ker (\psi ) = \{ 0 \} .$ Notice that
 $(q_1,i(q_1),q_3,0)\in {\cal B}_3$. Using $f=(id_F, \psi )$, Theorem 32 in
 \cite{ludkmicrta19}, Remark 11 and Corollary 20, we deduce  that there exists $\phi \in Hom({\cal
 B}_1,{\cal F})$ such that $\phi ({\cal G})= \{ (0,0,0,0) \} $ and
 $\phi (q)=(0,r,q_3,r)$.

\par {\bf Theorem 23.} {\it Let ${\cal B}={\cal
B}(A,E,F,X_F,i,p)$ and ${\cal R}={\cal B}(A',E',F,Y_F,i',p')$ be
completely $(F,r)$-regular microbundles and let ${\cal R}$ be also
$(F,r)$-arcwise connected. Let $\alpha $ and $\delta $ be closed
$M$-networks on ${\cal B}$. Then $\tau _{\alpha }\subset \tau
_{\delta }$ for $M$-set-open topologies on $Hom ({\cal B},{\cal R})$
if and only if for each $C\in \alpha $ there exist $C_1$,...,$C_n$
in $\delta $, $n\in {\bf N}$, such that $C\subset C_1 \cup ... \cup
C_n$.}
\par {\bf Proof.} Assume that for each $C\in \alpha $ there exist $C_1$,...,$C_n$ in
$\delta $, $n\in {\bf N}$, such that $C\subset C_1 \cup ... \cup
C_n$. Then from Definitions 1, 16, 18 and Remarks 11, 19 it follows
that $\tau _{\alpha }\subset \tau _{\delta }$. \par Vice versa,
assume that $\tau _{\alpha }\subset \tau _{\delta }$. By the
conditions of this theorem the microbundle ${\cal R}$ is
$(F,r)$-arcwise connected. Let $q\in {\cal R}$, $s\in {\cal R}$ with
$s_3=0$, $s_4=0$, $q\ne s$, $(q_1,q_2)\ne (s_1,s_2)$. Then there
exists $\mu =\mu _{q,s} \in Hom ({\cal F},{\cal R})$ such that
$s=\mu (0,0,0,0)$ and $q=\mu (0,r,q_3,r)$ according to Definition
21. We take any $C\in \alpha $ and an open submicrobundle ${\cal
G}={\cal B}(A_0,E_0,F,X_{0,F},i_0,p_0)$ in ${\cal R}$ such that
$s\in {\cal G}$ and $q'=(q_1,q_2,q_3',q_4)\notin {\cal G}$ for each
$q_3'\in F$, where $i_0=i'|_{A_0}$, $p_0=p'|_{E_0}$. Therefore,
${\sf H}(C,{\cal G})$ is a neighborhood of $\psi \in Hom_{\delta
}({\cal B},{\cal R})$ such that $\psi ({\cal B})=(s_1,s_2,0,0)$ with
$\psi ^1_1(A)= \{ s_1 \} $, $\psi ^2_2(E)= \{ s_2 \} $, $\psi
^3_3(F)= \{ 0 \} $, $\psi ^4_4(X_F)=\{ 0 \} $, where $s_1\in A$,
$s_2\in E$, $i(s_1)=s_2$.
\par From the condition $\tau _{\alpha }\subset \tau _{\delta
}$ it follows that there exist $C_1$,...,$C_n$ in $\delta $, $n\in
{\bf N}$, open subsets $D_1$,...,$D_n$ in ${\cal R}$, such that
$P\subset {\sf H}(C, {\cal G})$, $\psi \in P$, $P:= {\sf
H}(C_1,D_1)\cap ... \cap {\sf H}(C_n,D_n)$. Assume that there exists
$g\in C$, $g\notin {\cal J}$, where ${\cal J}:=C_1 \cup ... \cup
C_n$. The microbundle ${\cal B}$ is completely $(F,r)$-regular,
consequently, there exists $\phi \in Hom ({\cal B},{\cal F})$ such
that $\phi ({\cal J})= \{ (0,0,0,0) \} $ and $\phi (g)=(0,r,g_3,r)$.
Then $\mu \circ \phi \in P$, $\mu \circ \phi \notin {\sf H}(C,{\cal
G})$ providing a contradiction. Thus $C\subset {\cal J}$.

\par {\bf Examples 24.(1).} Let ${\cal B}$ be the
$T_1$-microbundle and let $p$ be a closed $M$-network consisting of
all closed subsets $C$ of ${\cal B}$ of the form $C=C_1\cup ... \cup
C_n$, where $C_j\in \alpha _0$ for each $j$, $n\in {\bf N}$, with
$\beta = (1,2,3,4)$ (see Definitions 1, 18, Remarks 2, 19). Then
$Hom_{p, \beta , \gamma } ({\cal B}, {\cal R})$ is supplied with the
topology $\tau _{p,\beta ,\gamma }$, which can be naturally called
$M$-point $\gamma $-open topology and denoted by $\tau _{p,\gamma
}$, or shortly $M$-point-open topology for $\gamma = \{ 1, 2, 3, 4
\} $ and denoted by $Hom_p({\cal B}, {\cal R})$ and $\tau _p$
respectively.
\par $\bf (2).$ If $c$ is a $\beta $-closed $M$-network consisting
of all $\beta $-compact subsets of ${\cal B}$, then $\tau _{co,
\beta , \gamma }$ is the $M$-$\beta $-compact-$\gamma $-open
topology on $Hom_{co, \beta , \gamma } ({\cal B}, {\cal R})$, or
shortly compact-open topology for $\beta = \gamma = \{ 1, 2, 3, 4 \}
$.
\par $\bf (3).$ Let ${\cal B}^{\tau }={\cal
B}(A,E,F,X_F,i,p)$ be the microbundle with the topology $\tau _A$ on
$A$, $\tau _E$ on $E$, $\tau _F$ on $F$, $\tau _X$ on $X$. Assume
that $\xi _A$ is a topology on $A$ finer than $\tau _A$, $\tau
_A\subset \xi _A$. Then there exists a coarsest topology $\xi _E$ on
the set $E$ such that $\tau _E\subset \xi _E$, $i: (A, \xi _A)\to
(E, \xi _E)$ is continuous, for each $b\in A$ there exists an open
neighborhood $U_{b,\xi }\in \xi _A$ of $b$ and an open neighborhood
$V_{b,\xi }\in \xi _E$ of $i(b)$ and a homeomorphism $h_{V_{b,\xi
}}: V_{b,\xi }\to U_{b,\xi }\times X$, where $\xi _F=\tau _F$, $\xi
_X=\tau _X$. Then there exists the microbundle ${\cal B}^{\xi }$
with the base space $(A, \xi _A)$, the total space $(E, \xi _E)$,
for the topological ring $(F, \tau _F)$ and the topological left
$F$-module $(X_F, \tau _X)$, such that $id=\pi : {\cal B}^{\xi }\to
{\cal B}^{\tau }$ is the continuous homomorphism, where $\pi
^1_1=id_A$, $\pi ^2_2=id_E$, $\pi ^3_3=id_F$, $\pi ^4_4=id_X$.
\par $\bf (4).$ For the microbundles ${\cal B}={\cal
B}(A,E,F,X_F,i,p)$ and ${\cal R}={\cal B}(A',E',F,X'_F,i',p')$
assume that $\xi _A=d_A$ is the discrete topology on $A$. Then
$Hom({\cal B}^d,{\cal R})$ is contained in $\Omega _{{\cal B},{\cal
R}}:= {A'}^A\times {E'}^E\times Hom (F,F) \times Hom (X,X')$ and the
latter is contained in $\Upsilon _{{\cal B},{\cal R}}:= {A'}^A\times
{E'}^E\times F^F \times {X'}^X$, where as usually $S^{\Lambda
}=\prod _{\lambda \in \Lambda } S_{\lambda }$, $S_{\lambda }=S$ for
each $\lambda \in \Lambda $, $S$ is a topological space, $S^{\Lambda
}$ is supplied with the Tychonoff product topology. Therefore, the
Tychonoff product topology on $\Upsilon _{{\cal B},{\cal R}}$
induces a topology on $Hom({\cal B}^d,{\cal R})$ which will be
denoted by $t_y=t_y({\cal B},{\cal R})$.
\par $\bf (5)$. The finest (i.e. largest ) $M$-set-$(\beta , \gamma
)$-open topology on $Hom ({\cal B},{\cal R})$ will be denoted by
$\tau _{l, \beta , \gamma }$, hence $\tau _{\alpha , \beta , \gamma
}\subset \tau _{l, \beta , \gamma }$ (see also Definition 18).

\par {\bf Corollary 25.} {\it Assume that ${\cal B}$ and ${\cal R}$
are $T_1$ microbundles. Then $Hom_p({\cal B},{\cal R})$ is dense in
$Hom_{t_y}({\cal B}^d,{\cal R})$.}
\par {\bf Proof.} This follows from Examples 24 $(3)$, $(4)$ and
Definition 18, since $C_p(A,A')$ is dense in ${A'}^A$, where
$C_p(A,A')$ denotes the space of all continuous mappings from $A$
into $A'$ supplied with the point-open topology.

\par {\bf Corollary 26.} {\it Assume that the conditions of Theorem 23 are
satisfied and $\alpha $ is a closed $M$-network on ${\cal B}$. Then
$\tau _{p,\gamma }\subset \tau _{\alpha , \gamma }$.}
\par {\bf Proof.} This follows from Theorem 23 and Example 24 $(1)$.

\par {\bf Definition 27.} If $\alpha $ is a $\beta $-closed
$M$-network on a microbundle ${\cal B}={\cal B}(A,E,F,X_F,i,p)$,
such that for each ${\cal G}\in \alpha $, ${\cal H}\subset _M{\cal
G}$ and ${\cal H}$ is $\beta $-closed in ${\cal B}$ implies ${\cal
H}\in \alpha $, then $\alpha $ is called a hereditarily $\beta
$-closed $M$-network (see also Definition 16, where ${\cal B}, {\cal
G}, {\cal H}$ are over the same topological ring $F$).
\par For a microbundle ${\cal R}={\cal
B}(A',E',F,X'_F,i',p')$ a family $\eta $ of $\gamma $-open subsets
will be called a $\gamma $-open $M$-base, if for each $\gamma $-open
subset ${\cal K}$ in ${\cal R}$, ${\cal K}\subset {\cal R}$, there
exists a family $ \{ {\cal A}_{\lambda }\in \eta : \lambda \in
\Lambda \} $, where $\Lambda $ is a set, such that ${\cal
K}=\bigcup_{\lambda \in \Lambda }{\cal A}_{\lambda }$. A family
$\zeta $ of $\gamma $-open subsets in ${\cal R}$ will be called a
$\gamma $-open $M$-subbase, if a family $ \{ {\cal A}_1 \cap  ...
\cap  {\cal A}_n\ne \emptyset : n\in {\bf N}, ~ \forall j\in \{
1,..., n  \} , {\cal A}_j\in \zeta \} $ is a $\gamma $-open $M$-base
for ${\cal R}$.

\par {\bf Theorem 28.} {\it Let ${\cal B}={\cal B}(A,E,F,X,i,p)$ and
${\cal R}={\cal B}(A',E',F,X',i',p')$ be $T_2$ microbundles, let
also the topological ring $F$ be commutative and associative, where
$X$ and $X'$ are topological left $F$-modules. Then $(Hom_F({\cal
B},{\cal R}), \tau _{co})$ is the $T_2$ microbundle; $(Hom({\cal
B},{\cal R}), \tau _{co})$ is the $T_2$ topological space.}
\par {\bf Proof.} Note that for $T_2$ topological spaces $J$, $K$
the topological space ${\sf C}(J,K)$ of all continuous maps from $J$
into $K$ supplied with the compact-open topology $\tau _{co}$ is
$T_2$ (see Sect. 44.I, Remark 1 in \cite{kurb}). Therefore,
$Hom_F({\cal B},{\cal R})$ and $Hom({\cal B},{\cal R})$ relative to
the $ \tau _{co}$ topology are $T_2$ as the topological spaces. \par
It remains to prove that $Hom_F({\cal B},{\cal R})$ can be supplied
with a microbundle structure. We recall that for each $f$, $g$ in
$Hom_F(X,X')$, $x$ and $y$ in $X$, $a$ and $b$ in $F$ the following
identities
\par $f(ax+by)=af(x)+bf(y)$,
\par $(f+g)(x)=f(x)+g(x)$,
\par $(bf)(x)=b(f(x))$,
\par $(a+b)f(x)=af(x)+bf(x)$,
\par $(ab)f(x)=a(bf(x))$,
\par $b(f(x))=bf(x)$,
\par $af(bx)=b(af(x))$
\\ imply that $Hom_F(X,X')$ is the left $F$-module,
since $X$ and $X'$ are the left $F$-modules, and $F$ is the
commutative associative ring.
\par By virtue of Theorem 1.1.7 in \cite{coyntb} $({\sf C}(X,X'), \tau
_{co})$ has a structure of and additive commutative topological
group. From the inclusion $Hom_F(X,X')\subset {\sf C}(X,X')$ and the
proof above it follows that $(Hom_F(X,X'), \tau _{co})$ also has the
structure of the additive commutative topological group.
\par For the topological left $F$-module $X'$ multiplication
$F\times X'\to X'$ is (jointly) continuous such that for each open
neighborhood $U'$ of $0$ in $X'$ there exist an open neighborhood
$V'$ of $0$ in $X'$ and an open neighborhood $S$ of $0$ in $F$ such
that $SV'\subset U'$. For a base $B_d(F)$ of open neighborhoods of
$d$ in $F$ and $e\in F$, $B_d(F)+e$ is a base of open neighborhoods
of $d+e$ in $F$. For a base $B_y(X')$ of open neighborhoods of $y$
in $X'$ and $z$ in $X'$, $B_y(X')+z$ is a base of open neighborhoods
of $y+z$ in $X'$. Notice that for each $b_0\in F$ and for each open
neighborhood $U'_1$ of $0$ in $X'$ there exists an open neighborhood
$V'_{b_0}$ of $0$ in $X'$ such that $b_0V'_{b_0}\subset U'_1.$ On
the other hand, for each $y_0\in X'$ and for each open neighborhood
$U'_2$ of $0$ in $X'$ there exists an open neighborhood $S_{y_0}$ of
$0$ in $F$ such that $S_{y_0}y_0\subset U'_2$. Then
\par $(28.1)$ $(b_0+S_1)(y_0+V'_1)=b_0y_0+S_1V'_1+b_0V'_1+S_1y_0\subset b_0y_0+U'+U'_1+U'_2$
for $S_1=S\cap S_{y_0}$, $V_1=V'_{b_0}\cap V'$. \par Let $Q$ be any
open  neighborhood of $0$ in $X'$ and $W$ be any compact subset in
$F$. For each $b\in W$ there exists an open neighborhood $J_b$ of
$b$ in $F$ and an open  neighborhood $P_b$ of $0$ in $X'$ such that
$J_bP_b\subset Q$ by $(28.1)$ with $U'$, $U'_1$, $U'_2$ such that
$U'+U'_1+U'_2\subset Q$, since $X'$ is the topological left
$F$-module. The open covering $ \{ J_b: b\in W \} $ has a finite
subcovering $ \{ J_{b_j}: j=1,...,n \} $, $n\in {\bf N}$,
$\bigcup_{j=1}^nJ_{b_j}\supset W$. Then $P=\bigcap_{j=1}^n P_{b_j}$
is the open neighborhood of $0$ in $X'$ and one gets that
$(\bigcup_{j=1}^nJ_{b_j})P\subset Q$, hence $WP\subset Q$.
\par For each $f$, $g$ in $Hom_F(X,X')$ such that $f-g\in {\sf
H}_4(B,P)$ there exist open in $(Hom_F(X,X'), \tau _{co})$ subsets
$\Omega _f$, $\Omega _g$ such that $\Omega _f-\Omega _g\subset {\sf
H}_4(B,P)$, consequently, $W(\Omega _f-\Omega _g)\subset {\sf
H}_4(B,WP)\subset {\sf H}_4(B,Q)$, since $(f-g)(WB)=W(f-g)(B)$ by
the left $F$-linearity of $f$ and $g$. By virtue of Theorem 31 this
implies that $(Hom_F(X,X'), \tau _{co})$ is the topological left
$F$-module. Similarly it is proved that if $Y$ is a topological left
$F$-module, then $({\sf C}(A,Y), \tau _{co})$ is the topological
left $F$-module, where $(b_1y_1+b_2y_2)(a)=b_1y_1(a)+b_2y_2(a)$ for
each $b_1$, $b_2$ in $F$, $y_1$, $y_2$ in ${\sf C}(A,Y)$, $a\in A$.
In particular, there $Y=(Hom_F(X,X'), \tau _{co})$ can be taken.
\par Next we consider $Hom_F(E,E'):=
\{ \pi ^2_2: \pi =(\pi ^1_1, \pi ^2_2, id_F, \pi ^4_4)\in
Hom_F({\cal B}, {\cal R}) \} $. For each $f\in {\sf C}(A,A')$ we put
$\tilde{i}(f)=i'\circ f\circ p$, consequently, $\tilde{i}(f)\in {\sf
C}(E,E')$. Then for each $g\in {\sf C}(E,E')$ let
$\tilde{p}(g)=p'\circ g\circ i$, hence $\tilde{p}(g)\in {\sf
C}(A,A')$.
\par If $f\in {\sf C}(A,A')$ and $g=\tilde{i}(f)$, then $g\in {\sf C}(E,E')$ as
the composition $g=i'\circ f \circ p$ of continuous maps. Let $C^2$
be compact in $E$ and $V'$ be open in $E'$ such that $g\in {\sf
M}(C^2,V'):= \{ t \in {\sf C}(E,E'): t(C^2)\subset V' \} $. Then
$p(C^2)=C^1$ is compact in $A$ by Theorem 3.1.10 in \cite{eng}. For
$U'$ open in $A'$ such that $V'\supset i'(U')$, the inclusion
$\tilde{i}({\sf M}(C^1,U'))\subset {\sf M}(C^2,V')$ is satisfied,
where ${\sf M}(C^1,U') = \{ s\in {\sf C}(A,A'): s(C^1)\subset U' \}
$, $f\in {\sf M}(C^1,U')$, $\tilde{i}(f)=g$. In view of Proposition
1.4.1 in \cite{eng} the map $\tilde{i}$ is continuous from $({\sf
C}(A,A'), \tau _{co})$ into $({\sf C}(E,E'), \tau _{co})$.
Similarly, if $g\in {\sf C}(E,E')$, $f=\tilde{p}(g)$, $C^1$ is a
compact subset in $A$, $U'$ is an open subset in $A'$, $V'$ is an
open subset in $E'$ such that $U'\supset p'(V')$, $g\in {\sf
M}(C^2,V')$, $C^2=i(C^1)$, then $\tilde{p}({\sf M}(C^2,V'))\subset
{\sf M}(C^1,U')$. Again by Theorem 3.1.10 and Proposition 1.4.1 in
\cite{eng} the map $\tilde{p}$ is continuous from $({\sf C}(E,E'),
\tau _{co})$ into $({\sf C}(A,A'), \tau _{co})$.
\par Notice that $\tilde{p}\circ \tilde{i}(f)=p'\circ i'\circ f
\circ p\circ i=f$ for each $f\in {\sf C}(A,A')$. If $f_1\ne f_2$ in
${\sf C}(A,A')$, then $f_1\circ p\ne f_2\circ p$ in ${\sf C}(E,A')$
and hence $i'\circ f_1\circ p\ne i'\circ f_2\circ p$ in ${\sf
C}(E,E')$, consequently, $\tilde{i}(f_1)\ne \tilde{i}(f_2)$.
\par From condition $(2.3)$ it follows that
\par $(28.2)$ $\pi ^1_1=p'\circ \pi ^2_2\circ i$ \\ for each $\pi \in
Hom_F({\cal B},{\cal R})$, since $p'\circ i'=id_{A'}$. Condition
$(2.4)$ implies that
\par $(28.3)$ $h'\circ \pi ^2_2\circ h^{-1}(a,ux+vy)=uh'\circ \pi ^2_2\circ
h^{-1}(a,x)+vh'\circ \pi ^2_2\circ h^{-1}(a,y)$ \\ for each $u$ and
$v$ in $F$, $x$ and $y$ in $X$, $a\in A$, since $\pi ^3_3=id_F$ in
the considered case, where $h'=h'_{V'}$ is for ${\cal R}$, $h=h_V$
is for ${\cal B}$. From condition $(2.5)$ we get that $\hat{\pi
}'_2\circ h'\circ \pi ^2_2=\pi ^4_4\circ \hat{\pi }_2\circ h$, where
$\hat{\pi }_2: U\times X\to X$, $\hat{\pi }'_2: U'\times X'\to X'$
by $(2.4)$, hence
\par $(28.4)$ $\hat{\pi }'_2\circ h'\circ \pi ^2_2\circ h^{-1}=\pi ^4_4\circ \hat{\pi
}_2$, \\ where $\pi ^4_4(ux+vy)=u\pi ^4_4(x)+v\pi ^4_4(y)$ for each
$x$ and $y$ in $X$, $u$ and $v$ in $F$, $\hat{\pi }_2(d,x)=x$ for
each $d\in U$ and $x\in X$.
\par Since ${\cal B}$ and ${\cal R}$ are $T_2$, then $X\setminus \{
0 \} $ is open in $X$, $X'\setminus \{ 0 \} $ is open in $X'$. Using
local homeomorphisms $h_V: V\to U\times X$ and $h'_{V'}: V'\to
U'\times X'$ for ${\cal B}$ and ${\cal R}$ respectively one gets
from $(I2.2)$, $(I2.3)$, that $i(A)$ is closed in $E$, $i'(A')$ is
closed in $E'$.
\par Let ${\cal B}_1={\cal B}(A,E_1,F,X',i_1,p_1)$ be a
$T_2$-microbundle. From $(28.2)$-$(28.4)$ we deduce that each $g\in
{\sf C}(A,A')$ induces $\pi \in Hom_F({\cal B}_1,{\cal R})$ such
that $\pi ^1_1=g$, $\pi ^3_3=id_F$, $\pi ^4_4=id_{X'}$, $\pi =(\pi
^1_1, \pi ^2_2, \pi ^3_3, \pi ^4_4)$, where $id_{X'}(y)=y$ for each
$y\in X'$, since $ \{ V_{1,i_1(b)}: b\in A \} $ is an open covering
of $E_1$, $ \{ V'_{i'(b')}: b'\in A' \} $ is an open covering of
$E'$, where $h_{1,V_{1,i_1(b)}}: V_{1,i_1(b)}\to U_b\times X'$ is a
homeomorphism for ${\cal B}_1$, $U_b$ is an open neighborhood of $b$
in $A$, $V_{1,i_1(b)}$ is an open neighborhood of $i_1(b)$ in $E_1$
(see condition $(I2.3)$).
\par The latter and the proof above imply, that there exists a map
$\zeta : {\sf C}(A,A')\to Hom_F({\cal B}_1,{\cal R})$ such that
$\zeta (g)=(g, \pi ^2_2, id_F,id_{X'})$. Conditions $(2.1)$-$(2.5)$
imply, that if $\pi = (\pi ^1_1, \pi ^2_2, id_F,\pi ^4_4)\in
Hom_F({\cal B},{\cal R})$, then $\pi ^4_4\in Hom_F(X,X')$. Using
$\zeta $, compositions of homomorphisms, Theorem 32 in
\cite{ludkmicrta19} and Corollary 20 (with $G=F$ and $X_G=X'$,
$X_F=X$) we infer that
\par $(28.5)$ $Hom_F(X,X')= \{ (id_F, \pi ^4_4): \pi = (\pi
^1_1, \pi ^2_2, id_F,\pi ^4_4)\in Hom_F({\cal B},{\cal R}) \} $.
\par For each $q=(q^1,q^2,q^3,q^4)\in {\cal B}$ and
$s=(s^1,s^2,s^3,s^4)\in {\cal R}$ according to conditions
$(I2.2)$-$(I2.4)$
\par $(28.6)$ $p\circ i(q^1)=q^1$, $p'\circ i'(s^1)=s^1$,
\par $(28.7)$ $\iota _0(q^1)=(q^1,0)$, $\iota '_0(s^1)=(s^1,0)$,
\par $(28.8)$ $\hat{\pi }_1(q^1,x)=q^1$, $\hat{\pi }_2(q^1,x)=x$  for each $x\in
X$, \par $\hat{\pi }_1\circ h_{V_{i(q^1)}}(q^2)=q^1$, $\hat{\pi
}_2\circ h_{V_{i(q^1)}}(q^2)=q^4$, \par $\hat{\pi }'_1(s^1,y)=s^1$,
$\hat{\pi }'_2(s^1,y)=y$ for each $y\in X'$, \par $\hat{\pi
}'_1\circ h'_{V'_{i'(s^1)}}(s^2)=s^1$, $\hat{\pi }'_2\circ
h'_{V'_{i'(s^1)}}(s^2)=s^4$, \\ where $\hat{\pi }_1$, $\hat{\pi
}_2$, $\iota _0$, $h_{V_{i(q^1)}}$ correspond to ${\cal B}$,
$\hat{\pi }'_1$, $\hat{\pi }'_2$, $\iota '_0$, $h'_{V'_{i'(s^1)}}$
correspond to ${\cal R}$, where $h_{V_{i(q^1)}}: V_{i(q^1)}\to
U_{q^1}\times X$, $h'_{V'_{i'(s^1)}}: V'_{i'(s^1)}\to U'_{s^1}\times
X'$ are homeomorphisms, $V_{i(q^1)}$ is the open neighborhood of
$i(q^1)$ in $E$, $U_{q^1}$ is the open neighborhood of $q^1$ in $A$,
similarly $V'_{i'(s^1)}$ is the open neighborhood of $i'(s^1)$ in
$E'$, $U'_{s^1}$ is the open neighborhood of $s^1$ in $A'$.
\par Therefore, for each $B\subset {\cal B}$, $D\subset {\cal R}$,
$g\in {\sf H}(B,D)$ from $(28.4)$ it follows that
\par $\hat{\pi }'_2\circ h'_{V'}\circ g^2_2(B^2\cap V)=g^4_4\circ
\hat{\pi }_2\circ h_V(B^2\cap V)\subset D^4$, since $\hat{\pi
}_2\circ h_V(B^2\cap V)\subset B^4$, $\hat{\pi }'_2\circ
h'_{V'}(D^2\cap V')\subset D^4$. Notice also that $\hat{\pi }_1\circ
h_V(B^2\cap V)=B^1\cap U$, $\hat{\pi }'_1\circ h'_{V'}(D^2\cap
V')=D^1\cap U'$ by $(28.6)$-$(28.8)$, since $i(U)\subset V$ and
$i'(U')\subset V'$. In view of Theorems 3.4.1, 3.4.7 in \cite{eng}
and $(28.3)$-$(28.5)$ above there exists an isomorphism
\par $\xi _{U,U'} : (Hom_F(U\times X,U'\times X'), \tau _{co})\to
( \Upsilon _{U,U'}, \tau _{co})$ \\ with $\Upsilon _{U,U'}= \{
f=(f_1,f_2)\in {\sf C}(U,U'\times Hom_F(X,X')): f_1\in {\sf
C}(U,U'), f_2\in {\sf C}(U,Hom_F(X,X')) \} $ such that $\Upsilon
_{U,U'}\subset {\sf C}(U,U')\times {\sf C}(U,Hom_F(X,X'))$ for each
$U$ open in $A$ and $U'$ open in $A'$. For $B$ compact in ${\cal B}$
and $D$ open in ${\cal R}$ we put
\par $\tilde{V}= \{ \pi ^2_2: \pi ^2_2=\nabla _{b\in B^1, b'\in D^1}
(h'_{V'_{i'(b')}})^{-1}\circ g^2_2\circ h_{V_{i(b)}}|_{V_{i(b)}\cap
((g^2_2)^{-1}(U'_{b'}\times X'))}, g^2_2\in \xi ^{-1}_{A,A'}( \Xi
(B^1,D^1) ) \} $ with $\Xi (B^1,D^1) = \{ f=(f_1,f_2)\in \Upsilon
_{A,A'}: f_1\in {\sf M}(B^1,D^1), f_2\in {\sf C}(A,Hom_F(X,X')) \}
$, where $\nabla $ denotes a combination of maps (see Proposition
2.1.1 in \cite{eng}), since $(h'_{V'_{i'(b')}})^{-1}\circ g^2_2\circ
h_{V_{i(b)}}|_P= (h'_{V'_{i'(b'_1)}})^{-1}\circ g^2_2\circ
h_{V_{i(b_1)}}|_P$ with $P=P_{b,b',b_1,b'_1}(g^2_2)=V_{i(b)}\cap
((g^2_2)^{-1}(U'_{b'}\times X'))\cap V_{i(b_1)}\cap
((g^2_2)^{-1}(U'_{b'_1}\times X'))$ for each $b$ and $b_1$ in $A$,
$b'$ and $b'_1$ in $A'$ such that $P\ne \emptyset $. \par Therefore,
for $\pi\in {\sf H}_{F;{\cal B}, {\cal R}}(B,D)$ we get that $\pi
^1_1\in {\sf M}(B^1,D^1)=: \tilde{U}$, $\tilde{i}(\pi ^1_1)=\pi
^2_2\in \tilde{V}$ and there exists a homeomorphism
$\tilde{h}_{\tilde{V}}: \tilde{V}\to \tilde{U}\times {\sf
C}(A,Hom_F(X,X'))$. Then we put $\tilde{\hat{\pi
}}_1(\tilde{d},\tilde{x})=\tilde{d}$, $~ \tilde{\hat{\pi
}}_2(\tilde{d},\tilde{x})=\tilde{x}$, $ ~ \tilde{\iota
}_0(\tilde{d})=(\tilde{d},0)$, $~
\tilde{i}(\tilde{d})=\tilde{h}^{-1}_{\tilde{V}}(\tilde{d},0)$ for
each $\tilde{d}\in {\sf M}(B^1,D^1)$ and $\tilde{x}\in {\sf
C}(A,Hom_F(X,X'))$. Therefore, $\tilde{\hat \pi }_1\circ
\tilde{\iota }_0|_{\tilde{U}}=\tilde{p}|_{\tilde{V}}\circ
\tilde{i}|_{\tilde{U}}$. Thus $(Hom_F({\cal B},{\cal R}), \tau
_{co})$ is isomorphic with the microbundle ${\cal B}={\cal
B}(A_1,E_1,F,X_1,\tilde{i},\tilde{p})$, where $A_1={\sf C}(A,A')$,
$E_1=Hom_F(E,E')$, $X_1={\sf C}(A,Hom_F(X,X'))$ in the compact-open
topology on them.

\par {\bf Definition 29.} Let ${\cal B}={\cal B}(A,E,F,X,i,p)$ and
${\cal R}={\cal B}(A',E',F',i',p')$ be microbundles and $P\subset
Hom ({\cal B},{\cal R})$, $P= \{ \pi _{\lambda }\in Hom({\cal
B},{\cal R}): \lambda \in \Lambda \} $, where $\Lambda $ is a set.
We put $\mu _P\in Hom({\cal B},{\cal B}_1)$ with ${\cal B}_1={\cal
B}(A_1,E_1,F_1,X_1,i_1,p_1)$ to be $\mu _P((b,e,f,x))= \{ (\pi
^1_{1,\lambda }(b), \pi ^2_{2,\lambda }(e),\pi ^3_{3,\lambda
}(f),\pi ^4_{4,\lambda }(x)): \lambda \in \Lambda \} $ for each
$(b,e,f,x)\in {\cal B}$, where $A_1={A'}^{\Lambda }$,
$E_1={E'}^{\Lambda }$, $F_1={F'}^{\Lambda }$, $X_1={X'}^{\Lambda }$,
where $Y^{\Lambda }$ is supplied with the Tychonoff product topology
for a topological space $Y$. It will be said that $P$ separates
elements of ${\cal B}$ from closed sets (or submicrobundles) if for
each closed subset $C$ (or submicrobundle) in ${\cal B}$ and $q\in
{\cal B}\setminus C$ there exists $\pi _{\lambda }\in P$ such that
$\pi _{\lambda }(q)\notin cl_{\cal R}\pi _{\lambda }(C)$, where
$cl_{\cal R}\pi _{\lambda }(C)$ denotes the closure of $\pi
_{\lambda }(C)$ in ${\cal R}$. Then, $P$ separates elements of
${\cal B}$ if for each $q\ne s$ in ${\cal B}$ there exists $\pi
_{\lambda }\in P$ such that $\pi _{\lambda }(q)\ne \pi _{\lambda
}(s)$.

\par {\bf Theorem 30.} {\it Assume that $P\subset Hom({\cal B},{\cal
R})$, ${\cal B}$ and ${\cal R}$ are $T_1$ microbundles, $P$
separates elements of ${\cal B}$ from closed sets (see Definition
29). Then $\mu _P: {\cal B}\to {\cal B}_1$ is an embedding.}
\par {\bf Proof.} In view of Theorems 27 and 29 in \cite{ludkmicrta19}
${\cal R}^{\Lambda }$ is isomorphic with ${\cal B}_1$, where ${\cal
R}^{\Lambda }=\prod_{\lambda \in \lambda } {\cal R}_{\lambda }$,
where ${\cal R}_{\lambda }={\cal R}$ for each $\lambda \in \Lambda
$. Consider any open in ${\cal B}$ subset $W$ and $q\in W$.
Therefore, there exists $\pi _{\lambda }\in P$ such that $\pi
_{\lambda }(q)\notin cl_{\cal R}\pi _{\lambda }({\cal B}\setminus
W)$ by the conditions of this theorem. Let $S={\cal R}\setminus
cl_{\cal R}\pi _{\lambda }({\cal R}\setminus W)$ and let $T=\rho
_{\lambda }^{-1}(S)$, where $\rho _{\lambda }\in Hom ({\cal R},{\cal
B}_1)$ is induced by the isomorphism of ${\cal R}^{\Lambda }$ with
${\cal B}_1$ and the projection homomorphism from ${\cal R}^{\Lambda
}$ onto ${\cal R}_{\lambda }$. Take any $s\in {\cal B}$ such that
$\mu _P(s)\in T$.
\par Therefore, $s\in \pi _{\lambda }^{-1}(S)$ and $\pi _{\lambda }^{-1}(S)={\cal B}
\setminus \pi _{\lambda }^{-1}(cl_{\cal R}\pi _{\lambda }({\cal
B}\setminus W))\subset W$. This implies that $T\cap \mu _P({\cal
B})$ is an open neighborhood of $\mu _P(q)$ and $T\cap \mu _P({\cal
B})\subset \mu _P(W)$, consequently, $\mu _P(W)$ is open in $\mu
_P({\cal B})$.

\par {\bf Definition 31.} Let ${\cal B}_1$, ${\cal B}_2$, ${\cal
B}_3$ be microbundles, ${\cal B}_j={\cal
B}(A_j,E_j,F_j,X_j,i_j,p_j)$ for each $j$. For each $\pi ^1_2\in
Hom({\cal B}_1,{\cal B}_2)$, $\pi ^2_3\in Hom ({\cal B}_2,{\cal
B}_3)$ we put $\Phi (\pi ^1_2, \pi ^2_3)=\pi ^2_3\circ \pi ^1_2$.
Then $(\pi ^2_3)_*(\pi ^1_2)=\pi ^2_3\circ \pi ^1_2$ and $(\pi
^1_2)^*(\pi ^2_3)=\pi ^2_3\circ \pi ^1_2$ are called induced
homomorphisms, where $(\pi ^2_3)_*: Hom ({\cal B}_1,{\cal B}_2)\to
Hom ({\cal B}_1,{\cal B}_3)$ and $(\pi ^1_2)^*: Hom ({\cal
B}_2,{\cal B}_3)\to Hom ({\cal B}_1,{\cal B}_3)$.

\par {\bf Theorem 32.} {\it Let ${\cal B}_1$, ${\cal B}_2$, ${\cal
B}_3$ be $T_2$ microbundles, let $\alpha _1$ be a compact
$M$-network on ${\cal B}_1$ and $\alpha _2$ be a closed $M$-network
forming a neighborhood base on ${\cal B}_2$. Then the map $\Phi :
Hom _{\alpha _1}({\cal B}_1,{\cal B}_2)\times Hom _{\alpha _2}({\cal
B}_2,{\cal B}_3)\to Hom _{\alpha _1}({\cal B}_1,{\cal B}_3)$ is
continuous.}
\par {\bf Proof.} Choose any $C\in \alpha _1$ and $D$ open in ${\cal
B}_3$, $\pi ^1_2\in Hom ({\cal B}_1,{\cal B}_2)$, $\pi ^2_3\in Hom
({\cal B}_2,{\cal B}_3)$ such that $\Phi (\pi ^1_2, \pi ^2_3)\in
{\sf H}_{{\cal B}_1,{\cal B}_3}(C,D)$. Therefore, $\pi ^2_3(\pi
^1_2(C))\subset D$ and $\pi ^1_2(C)\subset (\pi ^2_3)^{-1}(D)$.
Notice that for each $q\in \pi ^1_2(C)$ there exists $B_q\in \alpha
_2$ such that $q\in Int_{{\cal B}_2} (B_q)$ and $B_q\subset (\pi
^2_3)^{-1}(D)$, where $Int_{{\cal B}_2} (B_q)=\{ s\in {\cal B}_2:
s\in (Int_{A_2} (B_q^1), Int_{E_2} (B_q^2), Int_{F_2} (B_q^3),
Int_{X_2} (B_q^4)) \} ,$ $B_q=(B_q^1,B_q^2,B_q^3,B_q^4)$, $Int_PS$
denotes the interior of $S$ in a topological space $P$ for $S\subset
P$ (see Definitions 1, 18 and Remark 2). By virtue of Theorem 3.1.10
in \cite{eng} $\pi ^1_2(C)$ is compact,consequently, there exist
$q_1,...,q_n$ in $\pi ^1_2(C)$, $n\in {\bf N}$, such that $\pi
^1_2(C)\subset \bigcup_{i=1}^nInt_{{\cal B}_2}B_{q_i}=:W$. This
implies that $(\pi ^1_2,\pi ^2_3)\in \Upsilon $, where $\Upsilon :=
{\sf H}_{{\cal B}_1, {\cal B}_2}(C,W)\times (\bigcap_{i=1}^n{\sf
H}_{{\cal B}_2, {\cal B}_3}(B_{q_i},D))$ and $\Phi (\Upsilon
)\subset {\sf H}_{{\cal B}_1, {\cal B}_3}(C,D)$.

\par {\bf Theorem 33.} {\it Assume that $\pi ^2_3\in Hom({\cal B}_2,{\cal
B}_3)$. If $\alpha _1$ is a closed $M$-network on ${\cal B}_1$, then
$(\pi ^2_3)_*: Hom_{\alpha _1}({\cal B}_1,{\cal B}_2)\to Hom_{\alpha
_1}({\cal B}_1,{\cal B}_3)$ is continuous.}
\par {\bf Proof.} Choose any $C\in \alpha _1$ and $D$ open in ${\cal
B}_3$, $C\subset {\cal B}_1$, $D\subset {\cal B}_3$. Then $\pi
^1_2\in (\pi ^2_3)^{-1}_*({\sf H}_{{\cal B}_1,{\cal B}_3}(C,D))$ if
and only if $\pi ^2_3(\pi ^1_2(C))\subset D$. The latter inclusion
is equivalent to $\pi ^1_2\in {\sf H}_{{\cal B}_1,{\cal B}_2}(C,(\pi
^2_3)^{-1}(D))$. Thus $ (\pi ^2_3)^{-1}_*({\sf H}_{{\cal B}_1,{\cal
B}_3}(C,D))={\sf H}_{{\cal B}_1,{\cal B}_2}(C,(\pi ^2_3)^{-1}(D))$.

\par {\bf Theorem 34.} {\it Let $\pi ^2_3\in Hom({\cal B}_2,{\cal B}_3)$
be an embedding of a microbundle ${\cal B}_2$ into a microbundle
${\cal B}_3$. Let $\alpha _1$ be a closed $M$-network on ${\cal
B}_1$. Then $(\pi ^2_3)_*: Hom_{\alpha _1}({\cal B}_1,{\cal B}_2)\to
Hom_{\alpha _1}({\cal B}_1,{\cal B}_3)$ is an embedding.}
\par {\bf Proof.} Since $\pi ^2_3$ is an injection, then $(\pi
^2_3)_*$ is an injection, that is $(\pi ^2_3)_*(\xi )\ne (\pi
^2_3)_*(\eta )$ for each $\xi \ne \eta $ in $Hom({\cal B}_1,{\cal
B}_2)$, since there exists $s\in {\cal B}_1$ such that $\xi (s)\ne
\eta (s)$. Take any $C\in \alpha _1$ and $B$ open in ${\cal B}_2$,
$C\subset {\cal B}_1$, $B\subset {\cal B}_2$, such that ${\sf
H}_{{\cal B}_1,{\cal B}_2}(C,B)$ belongs to a subbase of the $\tau
_{\alpha _1}$ topology on $Hom ({\cal B}_1,{\cal B}_2)$. There
exists $D$ open in ${\cal B}_3$, $D\subset {\cal B}_3$, such that
$D\cap \pi ^2_3({\cal B}_2)=\pi ^2_3(B),$ since $\pi ^2_3$ is the
embedding of the microbundle ${\cal B}_2$ into the microbundle
${\cal B}_3$. Therefore, $(\pi ^2_3)_*^{-1}({\sf H}_{{\cal
B}_1,{\cal B}_3}(C,D))={\sf H}_{{\cal B}_1,{\cal B}_2}(C,(\pi
^2_3)^{-1}(D))={\sf H}_{{\cal B}_1,{\cal B}_2}(C,B)$, consequently,
$(\pi ^2_3)_*({\sf H}_{{\cal B}_1,{\cal B}_2}(C,B))= {\sf H}_{{\cal
B}_1,{\cal B}_3}(C,D)\cap (\pi ^2_3)_*(Hom_{\alpha _1}({\cal
B}_1,{\cal B}_2))$.

\par {\bf Definition 35.} A homomorphism of microbundles $\pi ^1_2:
{\cal B}_1\to {\cal B}_2$ is almost onto if $\pi ^1_2({\cal B}_1)$
is dense in ${\cal B}_2$, that is $\pi ^{1,1}_{1,2}(A_1)$ is dense
in $A_2$, $\pi ^{2,1}_{2,2}(E_1)$ is dense in $E_2$, $\pi
^{3,1}_{3,2}(F_1)$ is dense in $F_2$, $\pi ^{4,1}_{4,2}(X_1)$ is
dense in $X_2$.

\par {\bf Theorem 36.} {\it Assume that $\pi ^1_2\in Hom ({\cal B}_1,{\cal B}_2)$,
where ${\cal B}_1$ and ${\cal B}_2$ are $T_2$ microbundles. \par
$(i)$. If $\pi ^1_2$ is almost onto, then $(\pi ^1_2)^*$ is an
injection. \par $(ii)$. If $\alpha _1$ is a closed $M$-network on
${\cal B}_1$, $Hom ({\cal B}_1,{\cal B}_3)$ separates elements in
${\cal B}_1$ and $(\pi ^1_2)^*: Hom ({\cal B}_2,{\cal B}_3)\to Hom
_{\alpha _1}({\cal B}_1,{\cal B}_3)$ is almost onto, then $\pi ^1_2$
is an injection.}
\par {\bf Proof.} Take any $\xi $ and $\eta $ in $Hom ({\cal B}_2,{\cal
B}_3)$ such that $(\pi ^1_2)^*(\xi )=(\pi ^1_2)^*(\eta )$. For each
$s\in \pi ^1_2({\cal B}_1)$ there exists $q=q_s\in {\cal B}_1$ such
that $\pi ^1_2(q)=s$, consequently, $\xi (s)=(\pi ^1_2)^*(\xi
)(q)=(\pi ^1_2)^*(\eta )(q)=\eta (s)$. Therefore, $\xi =\eta $,
since $\pi ^1_2({\cal B}_1)$ is dense in ${\cal B}_2$, and since
${\cal B}_1$, ${\cal B}_2$ and ${\cal B}_3$ are the $T_2$
microbundles.
\par $(ii)$. Choose any $s_1\ne s_2$ in ${\cal B}_1$. Since $Hom ({\cal B}_1,{\cal
B}_3)$ separates elements in ${\cal B}_1$, then there exists $\xi
\in Hom ({\cal B}_1,{\cal B}_3)$ such that $\xi (s_1)\ne \xi (s_2)$.
Choose disjoint open neighborhoods $W_1$ and $W_2$ of $\xi (s_1)$
and $\xi (s_2)$ in ${\cal B}_3$. They exist, since ${\cal B}_3$ is
the $T_2$ microbundle. Therefore, ${\sf H}_{{\cal B}_1,{\cal B}_3} (
\{ s_1 \} , W_1)\cap {\sf H}_{{\cal B}_1,{\cal B}_3}(\{ s_2 \} ,
W_2)=:\Psi $ is an open neighborhood of $\xi $ in $Hom _{\alpha
_1}({\cal B}_1,{\cal B}_3)$. There exists $\eta $ in $Hom ({\cal
B}_2,{\cal B}_3)$ such  that $(\pi ^1_2)^*(\eta ) \in \Psi $, since
$(\pi ^1_2)^*$ is almost onto. This implies that $\eta (\pi
^1_2(s_1))\in W_1$ and $\eta (\pi ^1_2(s_2))\in W_2$, consequently,
$\pi ^1_2(s_1)\ne \pi ^1_2(s_2)$.

\par {\bf Conclusion 37.} The results of this article can be used for
further investigations of topological structure of microbundles,
families of their continuous homomorphisms, their applications in
geometry, abstract harmonic analysis, mathematical physics,
dynamical systems, other branches of mathematics and other sciences.
This also can be used for studies of operator bundles, group
bundles, sheaves, foliations of manifolds or microbundles.

\end{document}